\documentclass{amsart}
\usepackage{amssymb}

\usepackage{graphicx}
\usepackage{tikz}
\usepackage{float}
\usepackage{subfig}

\newtheorem{theorem}{Theorem}[subsection]
\theoremstyle{plain} \numberwithin{equation}{section}

\newtheorem{convention}[theorem]{Convention}
\newtheorem{corollary}[theorem]{Corollary}

\newtheorem{definition}[theorem]{Definition}
\newtheorem{example}[theorem]{Example}

\newtheorem{lemma}[theorem]{Lemma}

\newtheorem{problem}{Problem}
\newtheorem{proposition}[theorem]{Proposition}
\newtheorem{remark}[theorem]{Remark}

\begin{document}
\title{On the transport dimension of measures}
\author{Qinglan Xia, Anna Vershynina}
\address{University of California at Davis\\
Department of Mathematics\\
Davis,CA,95616}
\email{qlxia@math.ucdavis.edu}
\urladdr{http://math.ucdavis.edu/\symbol{126}qlxia}
\subjclass[2000]{Primary 49Q20, 51F99; Secondary 28E05, 90B06}
\keywords{optimal transport path, ramified optimal transportation, fractal
dimension of measures, transport dimension, irrigation dimension}
\thanks{This work is supported by an NSF grant DMS-0710714.}

\begin{abstract}
In this article, we define the transport dimension of probability measures
on $\mathbb{R}^m$ using ramified optimal transportation theory. We show that
the transport dimension of a probability measure is bounded above by the
Minkowski dimension and below by the Hausdorff dimension of the measure.
Moreover, we introduce a metric, called ``the dimensional distance", on the
space of probability measures on $\mathbb{R}^m$. This metric gives a
geometric meaning to the transport dimension: with respect to this metric,
we show that the transport dimension of a probability measure equals to the
distance from it to any finite atomic probability measure.
\end{abstract}

\maketitle

\section{Introduction}

The theory of ramified optimal transportation aims at finding an optimal
transport path between two given probability measures. One of the measures
is representing the source while the other is representing the target. A
transport path is typically in the form of a tree-shaped branching
structure. A natural question is: given a probability measure, is it
possible to transport it to a Dirac mass via a finite cost transport path?
The answer to this question crucially depends on the dimensional information
of the given measure.

In \cite{Solimini}, Devillanova and Solimini studied the irrigability
dimension of measures using optimal transportation theory. For any given
probability measure $\mu $ on $\mathbb{R}^{m}$, the irrigability dimension
of $\mu $ is defined by

\begin{equation*}
\dim _{I}\left( \mu \right) :=\inf_{0\leq \alpha <1}\left\{ \frac{1}{%
1-\alpha }:\text{ if }\mu \text{ is }\alpha -\text{irrigable}\right\} .
\end{equation*}%
The main theorem in \cite{Solimini} shows that%
\begin{equation}
\max \left\{ \dim _{H}\left( \mu \right) ,1\right\} \leq \dim _{I}\left( \mu
\right) \leq \max \left\{ 1,\dim _{M}\left( \mu \right) \right\}
\label{irrigation}
\end{equation}%
where $\dim _{H}\left( \mu \right) $ (or $\dim _{M}\left( \mu \right) $)
denotes the infimum of the Hausdorff dimension (or the Minkowski dimension,
respectively) of sets that $\mu $ is concentrated on. By definition, the
irrigability dimension $\dim _{I}\left( \mu \right) $ of a measure $\mu $
must be larger or equal to $1$ as the parameter $\alpha $ is in the range of
$[0,1)$.

In this article, we aim at removing the maximum constraint from (\ref%
{irrigation}) by using a different approach of ramified optimal
transportation (i.e. using optimal transport paths), and also allowing the
parameter $\alpha $ to be negative. This generalization allows us to
consider measures which have fractal dimensions (e.g. the Cantor measure)
less than $1 $. We introduce an analogous concept called ``the transport
dimension of $\mu $'', and show in theorem \ref{main theorem} that
\begin{equation*}
\dim _{H}\left( \mu \right) \leq \dim _{T}\left( \mu \right) \leq \dim
_{M}\left( \mu \right)
\end{equation*}%
with a slight modification of the definition of $\dim _{M}\left( \mu \right)
$. The major difference between $\dim _{I}\left( \mu \right) $ and $\dim
_{T}\left( \mu \right) $ is that $\dim _{T}\left( \mu \right) $ is allowed
to be less than 1. For instance, for the Cantor measure, we show that $\dim
_{T}\left( \mu \right) $ is $\frac{\ln 2}{\ln 3}$ which is exactly the
dimension of the Cantor set.

Moreover, we find a geometric meaning to the transport dimension by
introducing a new metric, called ``the dimensional distance'', on the space
of equivalent classes of probability measures. With respect to this metric,
the distance of any probability measure to a Dirac mass (or any atomic
probability measure) equals to the transport dimension of the measure. In
other words, the transport dimension of a measure quantitatively describes
how far the measure is from being an atomic measure.

We refer to the book \cite{book} as well as references there for a modern
account of optimal transportation with branching structures. In particular,
a partial list of most relevant works is listed in the reference: \cite%
{buttazzo},\cite{BCM},\cite{DH},\cite{gilbert},\cite{msm},\cite{paolini},%
\cite{white},\cite{xia1},\cite{xia2},\cite{xia3},\cite{xia4},\cite{xia5}, %
\cite{xld},\cite{zhangzhu} and of course \cite{Solimini}.

\textbf{Organization of the paper.} We first recall some basic concepts
about optimal transport paths in section 2 with some necessary
modifications. We show that for any $\alpha \in \left( -\infty ,1\right) $,
there exists a metric $d_{\alpha }$ on the space $\mathcal{A}\left( \mathbb{R%
}^{m}\right) $ of atomic probability measures. Then, we consider the metric
completion $\mathcal{P}_{\alpha }\left( \mathbb{R}^{m}\right) $ of $\mathcal{%
A}\left( \mathbb{R}^{m}\right) $ with respect to each $d_{\alpha }$. We show
that when a probability measure is concentrated on a set whose box dimension
is less than $\frac{1}{1-\alpha }$, then the measure belongs to $\mathcal{P}%
_{\alpha }$. After that, we begin to consider the transport dimension of
measures in section 3 with some comparison with other dimension of measures.
In particular, we show that the transport dimension of a measure is bounded
above by its Minkowski dimension and bounded below by its Hausdorff
dimension. We also show the transport dimension of the Cantor measure is $%
\frac{\ln 2}{\ln 3}$. In section 4, we consider the dimensional distance
between probability measures. The main result there says that the transport
dimension of a probability measure is given by the distance of the measure
to any atomic probability measure with respect to the dimensional distance.

\section{The $d_{\protect\alpha }$ metric on probability measures for $%
\protect\alpha \in \left( -\infty ,1\right) $}

\subsection{Transport paths between atomic measures}

We first recall some basic concepts about transport paths between measures
of equal total mass as studied in \cite{xia1}, with some necessary
modifications.

Recall that a (finite) atomic measure on $\mathbb{R}^{m}$ is in the form of
\begin{equation}
\mathbf{a=}\sum_{i=1}^{k}m_{i}\delta _{x_{i}}  \label{prob_meas_a}
\end{equation}%
with distinct points $x_{i}\in \mathbb{R}^{m}$, and positive real numbers $%
m_{i},$ where $\delta _{x}$ denotes the Dirac mass located at the point $x$.
The measure $\mathbf{a}$ is a probability measure if its mass $\left| \left|
\mathbf{a}\right| \right| :=\sum_{i=1}^{k}m_{i}=1$. Let $\mathcal{A}(\mathbb{%
R}^{m})$ be the space of all atomic probability measures on $\mathbb{R}^{m}$.

\begin{definition}
Given two atomic measures
\begin{equation}
\mathbf{a}=\sum_{i=1}^{k}m_{i}\delta _{x_{i}}\text{ and }\mathbf{b}%
=\sum_{j=1}^{\ell }n_{j}\delta _{y_{j}}  \label{prob_meas}
\end{equation}%
in $\mathbb{R}^{m}$of equal mass, a \textbf{transport path} from $\mathbf{a}$
to $\mathbf{b}$ is a weighted directed graph $G$ consisting of a vertex set $%
V(G)$, a directed edge set $E(G)$ and a weight function $w:E(G)\rightarrow
(0,+\infty )$ such that $\{x_{1},x_{2},...,x_{k}\}\cup
\{y_{1},y_{2},...,y_{\ell }\}\subset V(G)$ and for any vertex $v\in V(G)$,
\begin{equation}
\sum_{e\in E(G),e^{-}=v}w(e)=\sum_{e\in E(G),e^{+}=v}w(e)+\left\{
\begin{array}{c}
m_{i},\text{\ if }v=x_{i}\text{\ for some }i=1,...,k \\
-n_{j},\text{\ if }v=y_{j}\text{\ for some }j=1,...,\ell \\
0,\text{\ otherwise }%
\end{array}%
\right.  \label{path}
\end{equation}%
where $e^{-}$ and $e^{+}$ denotes the starting and ending endpoints of each
edge $e\in E(G)$.
\end{definition}

Note that the balanced equation (\ref{path}) simply means the conservation
of mass at each vertex. In terms of polyhedral chains, we simply have $%
\partial G=b-a$.

For any two atomic measures $\mathbf{a}$ and $\mathbf{b}$ on $\mathbb{R}^{m}$
of equal total mass, let
\begin{equation*}
Path(\mathbf{a},\mathbf{b})
\end{equation*}
be the space of all transport paths from $\mathbf{a}$ to $\mathbf{b}$. It is
easy to see to that $Path\left( \mathbf{a},\mathbf{b}\right) $ is always
nonempty.

\begin{definition}
For any $\alpha \leq 1$ and any transport path $G\in \text{Path}(\mathbf{a},%
\mathbf{b})$, we define
\begin{equation*}
\mathbf{M}_{\alpha }(G):=\sum_{e\in E(G)}w(e)^{\alpha }\text{length}(e).
\end{equation*}
\end{definition}

We now consider the following ramified optimal transport problem:

\begin{problem}
Given two atomic measures $\mathbf{a}$ and $\mathbf{b}$ of equal mass on $%
\mathbb{R}^{m}$ and $-\infty<\alpha<1$, find a minimizer of
\begin{equation*}
\mathbf{M}_{\alpha }(G)
\end{equation*}%
among all transport paths $G\in Path\left( \mathbf{a},\mathbf{b}\right)$.
\end{problem}

An $\mathbf{M}_{\alpha }$ minimizer in $Path(\mathbf{a},\mathbf{b})$ is
called an \textit{$\alpha -$optimal transport path} from $\mathbf{a}$ to $%
\mathbf{b}$.

Note that in \cite[Definition 2.2]{xia1}, we only allow that $0\leq \alpha
\leq 1$. But now, for the purpose of studying fractal dimension of measures
in this article, we also allow negative $\alpha $. Negative $\alpha $ also
corresponds to many transport problems in reality. For instance, the
``cost'' (i.e. risk here) for a baby to go back home from his preschool by
himself is much higher than the one for an adult. This is why babies need to
be picked up by adults. In this example, the cost of transporting a higher
density mass (i.e. an adult) is less than the cost of transporting a lower
density mass (i.e. a baby). Negative $\alpha $ is used for such phenomenon.
Another application of ramified optimal transportation with negative $\alpha$
may be found in \cite{xia6}.

\begin{definition}
For any $\alpha \leq 1$, we define
\begin{equation*}
d_{\alpha }\left( \mathbf{a},\mathbf{b}\right) =\inf \left\{ \mathbf{M}%
_{\alpha }\left( G\right) :G\in Path\left( \mathbf{a},\mathbf{b}\right)
\right\}
\end{equation*}%
for any $\mathbf{a},\mathbf{b\in }\mathcal{A}(\mathbb{R}^{m})$.
\end{definition}

\begin{remark}
If $\mathbf{\tilde{a}}$ and $\mathbf{\tilde{b}}$ are two atomic measures of
equal total mass $\Lambda $, and let $\mathbf{a}=\frac{1}{\Lambda }\mathbf{%
\tilde{a}}$ and $\mathbf{b}=\frac{1}{\Lambda }\mathbf{\tilde{b}}$ be the
normalization of $\mathbf{\tilde{a}}$ and $\mathbf{\tilde{b}}$. Then, we set
\begin{equation*}
d_{\alpha }\left( \mathbf{\tilde{a}},\mathbf{\tilde{b}}\right) =\Lambda
^{\alpha }d_{\alpha }\left( \mathbf{a},\mathbf{b}\right) .
\end{equation*}
\end{remark}

\subsection{Behavior of transport paths for a negative $\protect\alpha $}

The behavior of an $\alpha -$optimal transport path is quite different for a
negative $\alpha $ and a nonnegative $\alpha $. For instance, in %
\cite[proposition 2.1]{xia1}, we shown that for nonnegative $\alpha$, an $%
\alpha$-optimal transport path contains no cycles. On the hand, from the
following example \ref{example_1}, an $\alpha -$optimal transport path may
prefer to have a cycle when $\alpha $ is negative.

\begin{example}
\label{example_1} Let $\mathbf{a=\delta }_{0}$, $\mathbf{b=}0.8\delta
_{0}+0.2\delta _{1}\in \mathcal{A}(\mathbb{R})$, then we may construct two
transport paths from $\mathbf{a}$ to $\mathbf{b}$. One contains a cycle
while the other does not.

\begin{figure}[h]
\centering \includegraphics[width=\textwidth]{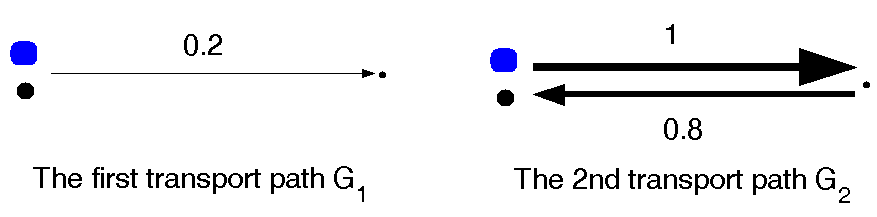}
%\label{negativefigure}
\caption{For $\protect\alpha<0$, an optimal transport path may contain a
cycle}
\label{cycle_example}
\end{figure}

The first transport path $G_{1}$ consists of only one directed edge from $0$
to $1$ with a weight $0.2$. The second transport path $G_{2}$ consists of
two directed edges: one from $0$ to $1$ with a weight $1$ and the other from
$1 $ to $0$ with a weight $0.8$. Then, when $\alpha <-\frac{1}{2}$, we have
\begin{equation*}
\frac{\mathbf{M}_{\alpha }\left( G_{2}\right) }{\mathbf{M}_{\alpha }\left(
G_{1}\right) }=\frac{1^{\alpha }+0.8^{\alpha }}{0.2^{\alpha }}\leq \frac{%
0.8^{\alpha }+ 0.8^{\alpha }}{0.2^{\alpha }}=2^{1+2\alpha }<1.
\end{equation*}%
Thus, $\mathbf{M}_{\alpha }\left( G_{2}\right) <\mathbf{M}_{\alpha }\left(
G_{1}\right) $. i.e. a path containing a cycle may have less $\mathbf{M}%
_{\alpha}$ cost.
\end{example}

A transport path may also contain another type of cycles, which are
unpleasant for our study.

\begin{figure}[h]
\centering \includegraphics[width=0.45\textwidth]{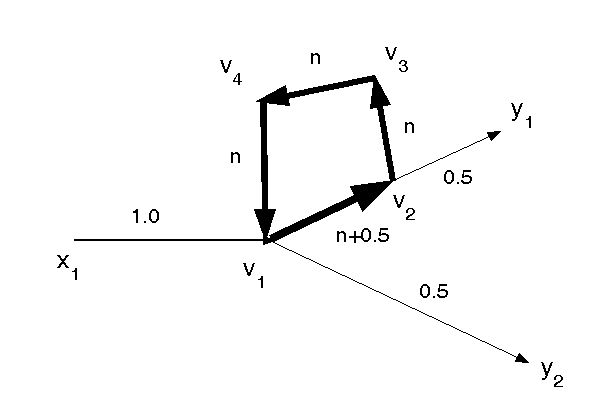}
\caption{An unpleasant cycle appeared in a transport path}
\label{directed_cyclic}
\end{figure}

\begin{example}
\label{example_2} From definition, it is possible for a transport path to
contain another type of directed cycles which are similar to the one
appearing in figure \ref{directed_cyclic}. Using this directed cycle, the
original weight `` $0.5$" of the edge in the example may either be
repeatedly counted or be combined with some mass from unknown sources (which
may have nothing to do with $\mathbf{a}$ or $\mathbf{b}$), and lead to a new
very high weight `` $n+0.5$" in the end. When $\alpha $ is negative and $n$
is large, $(n+0.5)^{\alpha}$ may become a very small positive number, much
less than $0.5^{\alpha}$. Thus, for negative $\alpha $, one may
significantly decrease the $\mathbf{M}_{\alpha }$ cost of a transport path
by repeatedly adding a directed cycle of this type with a very large
weighting constant $n$. Since the value $n^{\alpha }$ approaches zero as $n$
approaches to infinity when $\alpha <0$, allowing such kind of directed
cycles in the family $Path\left( \mathbf{a},\mathbf{b}\right) $ may
eventually lead to an unpleasant result: $d_{\alpha }\left( \mathbf{a},%
\mathbf{b}\right) =0$ for $\mathbf{a}\neq \mathbf{b}$.
\end{example}

The cycle in example \ref{example_1} is desirable while the cycle in example %
\ref{example_2} is not what we want. To overcome this conflict of interests,
we adopt the following convention on the definition of transport paths when $%
\alpha $ is negative.

\begin{convention}
\label{convention} For any vertex $v\in V\left( G\right) $ of a transport
path $G$, if there exists a list of vertices $\left\{ v_{1},v_{2,}\cdots
,v_{n}\right\} $ such that $v_{1}=v_{n}=v$ and $\left[ v_{i},v_{i+1}\right] $
is a directed edge in $E\left( G\right) $ with a positive edge length for
each $i=1,2,\cdots ,n-1$, then we will view $v_{n}$ as a \textbf{different}
copy of the point $v $. In other words, $v_{n}$ and $v$ are two different
vertices in $V\left( G\right) $ and thus the balance equation (\ref{path})
must be separately hold at each of them.
\end{convention}

Using this convention, one cannot ``combine" the weights of edges from the
vertex $v_n$ with weights of edges from the different vertex $v_1$. So the
weight of an edge cannot be repeatedly counted leading to a very large
number, and must be bounded above by the total mass of $\mathbf{a}$. As a
result, the convention implies two consequences: a transport path will no
longer contain a directed cycle in the type of example \ref{example_2}.
Moreover, we have a universal upper bound (i.e. the total mass of the
source) on the weight of each edge:
\begin{equation}  \label{w_e}
w\left( e\right) \leq \left| \left| \mathbf{a}\right| \right|
\end{equation}%
for each edge $e$ of a transport path $G$ from $\mathbf{a}$ to $\mathbf{b}$.
In particular, if $\mathbf{a}$ is a probability measure, then $w\left(
e\right) \leq 1$.

For nonnegative $\alpha $, since an $\alpha -$optimal transport path
contains no cycles (as shown in \cite[proposition 2.1]{xia1}), there is no
need to adopt this convention then.

Many properties about transport paths have been studied in \cite{xia1},\cite%
{xia2},\cite{xia7}, \cite{xia3},\cite{xia5} and \cite{book} when $0\leq
\alpha <1$. For instance, we showed in \cite[theorem 5.1]{xia1} that $\left(
\mathcal{A}(\mathbb{R}^{m}),d_{\alpha }\right) $ is a geodesic space, and it
is indeed a metric induced by a quasimetric (see \cite{xia5}). As the main
purpose of this article is studying dimensions of measures, we will leave
the study of properties of optimal transport paths in the situation of $%
\alpha <0$ to a later article. Currently, we only need to show that $%
d_{\alpha }$ is still a metric on $\mathcal{A}(\mathbb{R}^{m})$ when $\alpha
<0$.

\subsection{The $d_{\protect\alpha }$ metric when $\protect\alpha <0$}

We denote $S\left( p,r\right) $ (and $\bar{B}\left( p,r\right) $,
respectively) to be the sphere (and the closed ball, respectively) centered
at $p\in \mathbb{R}^{m}$ of radius $r>0$. Note that for any transport path $%
G $, the restriction of $G$ on any closed ball $\bar{B}\left( p,r_{0}\right)
$ gives another transport path $G|_{\bar{B}\left( p,r_{0}\right) }$.

\begin{lemma}
\label{estimate}Suppose $\mathbf{a}$ and $\mathbf{b}$ are two atomic
measures on $\mathbb{R}^{m}$ of equal total mass, and $G$ is a transport
path from $\mathbf{a}$ to $\mathbf{b}$. If the intersection of $G\cap
S\left( p,r\right) $ as sets is nonempty for almost all $r\in \left[ 0,r_{0}%
\right] $ for some $r_{0}>0$, then for any $\alpha <0$, we have
\begin{equation*}
\mathbf{M}_{\alpha }\left( G|_{\bar{B}\left( p,r_{0}\right) }\right) \geq
\Lambda ^{\alpha }r_{0},
\end{equation*}%
where $\Lambda $ is any upper bound of the weights of edges in $G|_{\bar{B}%
\left( p,r_{0}\right) }$.
\end{lemma}

\begin{proof}
For any $r>0$, let
\begin{equation*}
\mathbf{E}_{r}:=\left\{ e\in E\left( G\right) :e\cap S\left( p,r\right) \neq
\emptyset \right\}
\end{equation*}%
be the family of all edges of $G$ that intersects with the sphere $S\left(
p,r\right) $. By assumption, $\mathbf{E}_{r}$ is nonempty for almost all $%
r\in \left[ 0,r_{0}\right] $.

Let $L$ be any fixed ray with endpoint $p$, and let $P:\mathbb{R}%
^{m}\rightarrow L$ be the projection that maps any point $x\in \mathbb{R}%
^{m} $ to the point $P\left( x\right) \in L$ with $\left| P\left( x\right)
-p\right| =\left| x-p\right| $. Note that for any edge $e\in E\left(
G\right) $, under the projection $P$, the length of the segment $P\left(
e\right) $ is no more than the length of $e$. Therefore,
\begin{eqnarray*}
\mathbf{M}_{\alpha }(G|_{\bar{B}\left( p,r_{0}\right) }) &=&\sum_{e\in
E\left( G\mid B\left( p,r_{0}\right) \right) }\left[ w(e)\right] ^{\alpha }%
\text{length}(e) \\
&\geq &\sum_{e\in E\left( G\mid B\left( p,r_{0}\right) \right) }\left[ w(e)%
\right] ^{\alpha }\text{length}(P\left( e\right) ) \\
&=&\int_{0}^{r_{0}}\sum_{e\in E\left( G\mid B\left( p,r_{0}\right) \right) }
\left[ w(e)\right] ^{\alpha }\chi _{P\left( e\right) }\left( r\right) dr \\
&=&\int_{0}^{r_{0}}\sum_{e\in \mathbf{E}_{r}}\left[ w(e)\right] ^{\alpha }dr.
\end{eqnarray*}%
where in the second equality $\chi _{A}$ denotes the characteristic function
on a set $A$. So, we have shown that
\begin{equation}
\mathbf{M}_{\alpha }(G|_{\bar{B}\left( p,r_{0}\right) })\geq
\int_{0}^{r_{0}}\sum_{e\in \mathbf{E}_{r}}\left[ w(e)\right] ^{\alpha }dr.
\label{key_estimate}
\end{equation}%
Now, since $w\left( e\right) \leq \Lambda $ for any $e\in \mathbf{E}_{r}\neq
\emptyset $ and $\alpha <0$, (\ref{key_estimate}) yields
\begin{equation*}
\mathbf{M}_{\alpha }(G|_{\bar{B}\left( p,r_{0}\right) })\geq
\int_{0}^{r_{0}}\max_{e\in \mathbf{E}_{r}}\left[ w(e)\right] ^{\alpha
}dr\geq \int_{0}^{r_{0}}\Lambda ^{\alpha }dr=\Lambda ^{\alpha }r_{0}.
\end{equation*}
\end{proof}

\begin{corollary}
\label{negative_estimate}Suppose $\alpha \leq 0$ and $p\in \mathbb{R}^{m}$.
Then, for any atomic measure
\begin{equation*}
\mathbf{a=}\sum_{i=1}^{k}m_{i}\delta _{x_{i}}
\end{equation*}%
on $\mathbb{R}^{m}$ with mass $\left| \left| \mathbf{a}\right| \right|
:=\sum_{i=1}^{k}m_{i}>0$ and any $G\in Path\left( \mathbf{a,\left| \left|
\mathbf{a}\right| \right| \delta }_{p}\right) $ we have
\begin{equation*}
\mathbf{M}_{\alpha }\left( G\right) \geq \left| \left| \mathbf{a}\right|
\right| ^{\alpha }\left| p-x_{i}\right|
\end{equation*}%
for each $i=1,2,\cdots ,k$.
\end{corollary}

\begin{proof}
By (\ref{w_e}), $ w\left( e\right) \leq \left| \left| \mathbf{a}\right|
\right| $ for any edge $e\in E\left( G\right) $. Then, the result follows
from lemma \ref{estimate} by setting $r_{0}=\max_{1\leq i\leq k}\left\{
\left| p-x_{i}\right| \right\} $ and $\Lambda =$ $\left| \left| \mathbf{a}%
\right| \right| $.
\end{proof}

The inequality (\ref{key_estimate}) also gives a lower bound estimate for
positive $\alpha $.

\begin{corollary}
\label{positive_estimate}Suppose $0\leq \alpha <1$. For any $\mathbf{a\in }%
\mathcal{A}(\mathbb{R}^{m})$ in the form of (\ref{prob_meas_a}), $p\in
\mathbb{R}^{m}$ and $r_{0}>0$, we have
\begin{equation}
\left[ \sum_{d\left( p,x_{i}\right) >r_{0}}m_{i}\right] ^{\alpha }\leq \frac{%
d_{\alpha }(\mathbf{a},\delta _{p})}{r_{0}}.  \label{positive_estimate_eqn}
\end{equation}
\end{corollary}

\begin{proof}
Let $\lambda =\sum_{d\left( p,x_{i}\right) >r_{0}}m_{i}$. Let $G$ be any
transport path from $\mathbf{a}$ to $\delta _{p}$. Then, for any $0<r\leq
r_{0}$, we have
\begin{equation*}
\sum_{e\in \mathbf{E}_{r}}w(e)\geq \lambda .
\end{equation*}%
By the inequality (\ref{key_estimate}), since the function $f\left( x\right)
=x^{\alpha }$ is concave on $\left[ 0,1\right] $ when $0\leq \alpha <1$, we
have
\begin{eqnarray*}
\mathbf{M}_{\alpha }(G|_{\bar{B}\left( p,r_{0}\right) }) &\geq
&\int_{0}^{r_{0}}\sum_{e\in \mathbf{E}_{r}}\left[ w(e)\right] ^{\alpha }dr \\
&\geq &\int_{0}^{r_{0}}\left[ \sum_{e\in \mathbf{E}_{r}}w(e)\right] ^{\alpha
}dr\geq \int_{0}^{r_{0}}\lambda ^{\alpha }dr=\lambda ^{\alpha }r_{0}\text{.}
\end{eqnarray*}%
Therefore, we have (\ref{positive_estimate_eqn}).
\end{proof}

\begin{proposition}
$d_{\alpha }$ is a metric on $\mathcal{A}(\mathbb{R}^{m})$.
\end{proposition}

\begin{proof}
Obviously, we only need to consider the case $\alpha <0$. In this case, it
is clear that $d_{\alpha }$ is nonnegative and symmetric. Now, for any $%
\mathbf{a},\mathbf{b,c\in }\mathcal{A}(\mathbb{R}^{m})$ and any path $%
G_{1}\in Path\left( \mathbf{a},\mathbf{b}\right) $ and $G_{2}\in Path\left(
\mathbf{b},\mathbf{c}\right) $, let $G_{3}$ be the disjoint union of the
directed weighted graphs $G_{1}$ and $G_{2}$. That is,
\begin{equation*}
V\left( G_{3}\right) =\left( V\left( G_{1}\right) \setminus V(\mathbf{b}%
)\right) \coprod \left( V\left( G_{2}\right) \setminus V(\mathbf{b})\right)
\coprod V(\mathbf{b})
\end{equation*}
and $E\left( G_{3}\right) =E\left( G_{1}\right) \coprod E\left( G_{2}\right)$%
, where the symbol $\coprod $ denotes the disjoint union of sets, and $V(%
\mathbf{b})$ denotes the vertex set (i.e. the support) of the measure $%
\mathbf{b}$. Then, by using the convention \ref{convention}, $G_{3}$ is a
transport path from $\mathbf{a}$ to $\mathbf{c}$ with
\begin{equation*}
\mathbf{M}_{\alpha }\left( G_{3}\right) =\mathbf{M}_{\alpha }\left(
G_{1}\right) +\mathbf{M}_{\alpha }\left( G_{2}\right) .
\end{equation*}%
Thus, by taking infimum, we have the triangle inequality
\begin{equation*}
d_{\alpha }\left( \mathbf{a},\mathbf{c}\right) \leq d_{\alpha }\left(
\mathbf{a},\mathbf{b}\right) +d_{\alpha }\left( \mathbf{b},\mathbf{c}\right)
\text{.}
\end{equation*}%
Now, we only need to check that if $\mathbf{a}\neq \mathbf{b}$, then $%
d_{\alpha }\left( \mathbf{a},\mathbf{b}\right) >0$. We may assume that $%
\mathbf{a}$ and $\mathbf{b}$ are in the forms of (\ref{prob_meas}). If the
supports of $\mathbf{a}$ and $\mathbf{b}$ are different, i.e. $\left\{
y_{1},y_{2},\cdots ,y_{l}\right\} \neq \left\{ x_{1},x_{2},\cdots
,x_{k}\right\} $ as sets, we may assume that $y_{1}\notin \left\{
x_{1},x_{2},\cdots ,x_{k}\right\} $. In this case, we set $p=y_{1}$ and
\begin{equation*}
r_{0}:=\min \left\{ \left| p-x_{i}\right| :i=1,2,\cdots ,k\right\} >0\text{.}
\end{equation*}%
If $\left\{ x_{1},x_{2},\cdots ,x_{k}\right\} =\left\{ y_{1},y_{2},\cdots
,y_{l}\right\} $ as sets, then we may assume
\begin{equation*}
\mathbf{a}=\sum_{i=1}^{k}m_{i}\delta _{x_{i}}\text{ and }\mathbf{b}%
=\sum_{i=1}^{k}n_{i}\delta _{x_{i}}
\end{equation*}%
and $m_{1}\neq n_{1}$. In this case, we set $p=x_{1}$ and let
\begin{equation*}
r_{0}:=\min \left\{ \left| p-x_{i}\right| :m_{i}\neq n_{i}\text{, }%
i=2,\cdots ,k\right\} >0\text{.}
\end{equation*}%
In any of these two cases, for any transport path $G$ from $\mathbf{a}$ to $%
\mathbf{b}$, the intersection of $G$ with the sphere $S\left( p,r\right) $
is nonempty for any $0<r<r_{0}$. By lemma \ref{estimate},
\begin{equation*}
\mathbf{M}_{\alpha }\left( G\right) \geq r_{0}
\end{equation*}%
since $w\left( e\right) \leq 1$ for any edge $e\in E\left( G\right) $.
Therefore,
\begin{equation*}
d_{\alpha }\left( \mathbf{a},\mathbf{b}\right) \geq r_{0}>0\text{.}
\end{equation*}
\end{proof}

\subsubsection{The completion of $\mathcal{A}\left( \mathbb{R}^m\right) $
with respect to $d_{\protect\alpha }$}

\begin{definition}
For any $\alpha \in (-\infty ,1]$, let $\mathcal{P}_{\alpha }(\mathbb{R}%
^{m}) $ be the completion of the metric space $\mathcal{A}(\mathbb{R}^{m})$
with respect to the metric $d_{\alpha }$.
\end{definition}

We will simply write $\mathcal{P}_{\alpha }(\mathbb{R}^{m})$ as $\mathcal{P}%
_{\alpha }$.

Note that when $\alpha =1$, the metric $d_{1}$ is the usual Monge's distance
on $\mathcal{A}(\mathbb{R}^{m})$ and $\mathcal{P}_{1}$ is the space $%
\mathcal{P}$ of all probability measures on $\mathbb{R}^{m}$. Now, the
following lemma implies that each element in $\mathcal{P}_{\alpha }$ can be
viewed as a probability measure on $\mathbb{R}^{m}$.

\begin{lemma}
If $\beta <\alpha $, then $\mathcal{P}_{\beta }(\mathbb{R}^{m})\subseteq
\mathcal{P}_{\alpha }(\mathbb{R}^{m})$, and for all $\mu ,\nu $ in $\mathcal{%
P}_{\beta }(\mathbb{R}^{m})$ we have $d_{\beta }(\mu ,\nu )\geq d_{\alpha
}(\mu ,\nu )$.
\end{lemma}

\begin{proof}
Note that for any path $G$ between any two probability measures $\mathbf{a}$
and $\mathbf{b}$ in $\mathcal{A}(\mathbb{R}^{m})$ we have
\begin{equation*}
\mathbf{M}_{\alpha }(G)=\sum_{e\in E(G)}w(e)^{\alpha }\text{length}(e)\leq
\sum_{e\in E(G)}w(e)^{\beta }\text{length}(e)=\mathbf{M}_{\beta }(G).
\end{equation*}%
This is because $w(e)\leq 1$ for any $e\in E(G)$. Hence
\begin{eqnarray*}
d_{\alpha }(\mathbf{a},\mathbf{b}) &=&\inf \{\mathbf{M}_{\alpha }(G):G\in
\text{Path}(\mathbf{a},\mathbf{b})\} \\
&\leq &\inf \{\mathbf{M}_{\beta }(G):G\in \text{Path}(\mathbf{a},\mathbf{b}%
)\}=d_{\beta }(\mathbf{a},\mathbf{b}).
\end{eqnarray*}%
Therefore, any Cauchy sequence\ $\left\{ \mathbf{a}_{n}\right\} $ in $%
\mathcal{A}(\mathbb{R}^{m})$ with respect to the metric $d_{\beta }$ is also
a Cauchy sequence with respect to $d_{\alpha }$. Thus, $\mathcal{P}_{\beta }(%
\mathbb{R}^{m})\subseteq \mathcal{P}_{\alpha }(\mathbb{R}^{m})$ with
\begin{equation*}
d_{\alpha }(\mu ,\nu )\leq d_{\beta }(\mu ,\nu )
\end{equation*}%
for any $\mu ,\nu $ in $\mathcal{P}_{\beta }(\mathbb{R}^{m})$.
\end{proof}

As discussed in \cite{Solimini} the Minkowski-Bouligand dimension, also
known as box-counting dimension of a set $A$ is given by
\begin{eqnarray}
\dim _{box}(A) &=&\inf \{\beta >0:A\text{ can be covered by }C_{\beta
}\delta ^{-\beta }\text{ balls of }  \notag  \label{box} \\
&&\text{radius }\delta \text{ for all }\delta \text{ less than }1\text{ and
some constant }C_{\beta }\}.
\end{eqnarray}%
A positive Borel measure $\mu $ on $\mathbb{R}^{m}$ is said to be \textit{%
concentrated} on a Borel set $A$ if $\mu (\mathbb{R}^{m}\setminus A)=0$.

\begin{theorem}
Suppose $\mu $ is a probability measure concentrated on a subset $A$ of $%
\mathbb{R}^{m}$ with $\dim _{box}(A)<\frac{1}{1-\alpha }$ for some $0\leq
\alpha <1$, then $\mu \in \mathcal{P}_{\alpha }$.
\end{theorem}

\begin{proof}
Since $\dim _{box}(A)<\frac{1}{1-\alpha }$, we may fix a constant $\beta $
such that $\dim _{box}(A)<\beta <\frac{1}{1-\alpha }$. By (\ref{box}), for
every $n\in \mathbb{N}$, the set $A$ can be covered by balls $\left\{
B\left( x_{i}^{\left( n\right) },\frac{1}{2^{n}}\right) \right\}
_{i=1}^{N_{n}}$ of radius $\frac{1}{2^{n}}$ centering at $x_{i}^{\left(
n\right) }$, where the number $N_{n}$ of balls of radius $\frac{1}{2^{n}}$
is bounded above by $C_{\beta }\left( \frac{1}{2^{n}}\right) ^{-\beta }$ for
some constant $C_{\beta }$. Using Vitali's covering theorem (see for
instance, \cite[section 1.5.1]{evans}) we may assume that the balls $\left\{
B\left( x_{i}^{\left( n\right) },\frac{1}{5\cdot 2^{n}}\right) \right\}
_{i=1}^{N_{n}}$ are disjoint. For each ball $B\left( x_{i}^{\left( n\right)
},\frac{1}{2^{n}}\right) $, let
\begin{equation*}
\Phi _{i}^{\left( n\right) }=\left\{ x_{j}^{\left( n-1\right) }:B\left(
x_{i}^{\left( n\right) },\frac{1}{2^{n}}\right) \cap B\left( x_{j}^{\left(
n-1\right) },\frac{1}{2^{n-1}}\right) \neq \emptyset \right\}
\end{equation*}%
be the set of centers of all balls in the preceding cover that intersects
with $B\left( x_{i}^{\left( n\right) },\frac{1}{2^{n}}\right) $. Then one
can easily check that the cardinality $k_{i}^{\left( n\right) }$ of $\Phi
_{i}^{\left( n\right) }$ is less than $\left( \frac{17}{2}\right) ^{m}$.

Using these family of covers, we construct a Cauchy sequence $\left\{
\mathbf{a}_{n}\right\} $ for $\mu $ as follows. For each $n$, let
\begin{equation*}
\mathbf{a}_{n}=\sum_{i=1}^{N_{n}}m_{i}^{\left( n\right) }\delta
_{x_{i}^{(n)}},
\end{equation*}%
where the positive number $m_{i}^{\left( n\right) }=\mu (B\left(
x_{i}^{\left( n\right) },\frac{1}{2^{n}}\right) \setminus
\bigcup_{h=1}^{i-1}B\left( x_{h}^{\left( n\right) },\frac{1}{2^{n}}\right) )$%
. Then, for each $n>1$, we may construct a partition of $\mathbf{a}%
_{n}=\sum_{j=1}^{N_{n-1}}\mathbf{a}_{n,j}^{(n-1)}$ with respect to $\mathbf{a%
}_{n-1}$ by setting
\begin{equation*}
\mathbf{a}_{n,j}^{(n-1)}=\sum_{i=1}^{N_{n}}m_{i,j}^{\left( n\right) }\delta
_{x_{i}^{(n)}}
\end{equation*}%
where the number
\begin{equation*}
m_{i,j}^{\left( n\right) }=:\mu \left( \left( B\left( x_{i}^{\left( n\right)
},\frac{1}{2^{n}}\right) \setminus \bigcup_{h=1}^{i-1}B\left( x_{h}^{\left(
n\right) },\frac{1}{2^{n}}\right) \right) \cap B\left( x_{j}^{\left(
n-1\right) },\frac{1}{2^{n-1}}\right) \right)
\end{equation*}%
if $x_{j}^{(n-1)}\in \Phi _{i}^{\left( n\right) }$ and we set $%
m_{i,j}^{\left( n\right) }$ to be zero otherwise. Also, let $G_{n,i}^{n-1}$
be the path sending a mass $m_{i,j}^{\left( n\right) }$ from the center $%
x_{i}^{(n)}$ to every point $x_{j}^{(n-1)}\in \Phi _{i}^{\left( n\right) }$.
Then,

\begin{eqnarray*}
d_{\alpha }\left( \mathbf{a}_{n-1},\mathbf{a}_{n}\right) &\leq
&\sum_{i=1}^{N_{n}}\mathbf{M}_{\alpha }\left( G_{n,i}^{n-1}\right) \\
&\leq &\sum_{i=1}^{N_{n}}\sum_{j=1}^{N_{n-1}}\left( m_{i,j}^{\left( n\right)
}\right) ^{\alpha }|x_{i}^{(n)}-x_{j}^{(n-1)}| \\
&\leq &\sum_{i=1}^{N_{n}}k_{i}^{\left( n\right) }\left( m_{i}^{\left(
n\right) }\right) ^{\alpha }(\frac{1}{2^{n}}+\frac{1}{2^{n-1}}).
\end{eqnarray*}%
Now, \textit{by H\"{o}lder inequality, we have }%
\begin{eqnarray*}
&&d_{\alpha }\left( \mathbf{a}_{n-1},\mathbf{a}_{n}\right) \\
&\leq &\left( \frac{17}{2}\right) ^{m}\left( \frac{3}{2^{n}}\right) \left(
\sum_{i=1}^{N_{n}}m_{i}^{\left( n\right) }\right) ^{\alpha }\left(
\sum_{i=1}^{N_{n}}1\right) ^{1-\alpha }\text{, } \\
&\leq &\left( \frac{17}{2}\right) ^{m}(\frac{3}{2^{n}})\left( N_{n}\right)
^{1-\alpha }\text{, \textit{because }}\sum_{i=1}^{N_{n}}m_{i}^{(n)}=1 \\
&\leq &\frac{3}{2^{n}}\cdot \left( \frac{17}{2}\right) ^{m}\left( C_{\beta
}\left( \frac{1}{2^{n}}\right) ^{-\beta }\right) ^{1-\alpha }=3\cdot \left(
\frac{17}{2}\right) ^{m}\left( C_{\beta }\right) ^{1-\alpha }b^{n},
\end{eqnarray*}%
where $b=(\frac{1}{2})^{1-\beta (1-\alpha )}\in \left( 0,1\right) $, because
$\beta <\frac{1}{1-\alpha }$. As a result, $\left\{ \mathbf{a}_{n}\right\} $
is a Cauchy sequence representing $\mu $ with respect to the metric $%
d_{\alpha }$ and therefore $\mu \in \mathcal{P}_{\alpha }$.
\end{proof}

\begin{corollary}
\label{p=pa}Suppose $\mu $ is a probability measure on $\mathbb{R}^{m}$ with
a compact support. Then $\mu \in \mathcal{P}_{\alpha }$ for any $\alpha \in
\left( 1-\frac{1}{m},1\right) $.
\end{corollary}

\section{Dimension of measures}

In this section, we will study properties of measures that belong to a
special subset of $\mathcal{P}_{\alpha }$, and then define the transport
dimension of measures.

\subsection{$d_{\protect\alpha }-$\textbf{admissible Cauchy sequence}}

\begin{definition}
\label{admissible}Let $\{\mathbf{a}_{k}\}_{k=1}^{\infty }$be a sequence of
atomic measures of equal total mass in the form of
\begin{equation*}
\mathbf{a}_{k}=\sum_{i=1}^{N_{k}}m_{i}^{\left( k\right) }\delta
_{x_{i}^{\left( k\right) }}
\end{equation*}%
for each $k$.
\begin{figure}[h!]
\centering \includegraphics[width=0.4\textwidth]{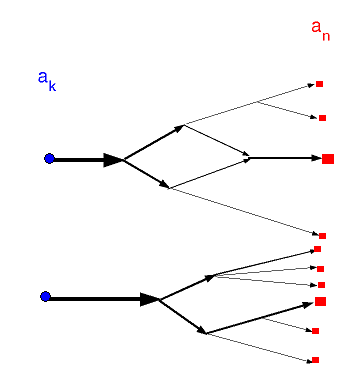}
\caption{An example of a transport path between $a_k$ and $a_n$.}
\label{admissible_path}
\end{figure}
We say that this sequence is a $d_{\alpha }-$\textbf{admissible Cauchy
sequence} if for any $\epsilon >0$, there exists an $N$ such that for all $%
n>k\geq N$ there exists a partition of
\begin{equation*}
\mathbf{a}_{n}=\sum_{i=1}^{N_{k}}\mathbf{a}_{n,i}^{(k)}
\end{equation*}%
with respect to $\mathbf{a}_{k}$ as sums of disjoint atomic measures and a
path (see figure \ref{admissible_path})
\begin{equation*}
G_{n,i}^{k}\in Path(m_{i}^{(k)}\delta _{x_{i}^{(k)}},\mathbf{a}_{n,i}^{(k)})
\end{equation*}%
for each $i=1,2,\cdots ,N_{k}$ such that
\begin{equation*}
\sum_{i=1}^{N_{k}}\mathbf{M}_{\alpha }\left( G_{n,i}^{k}\right) \leq
\epsilon \text{.}
\end{equation*}%
Also, we denote $G_{n}^{k}=\sum_{i=1}^{N_{k}}G_{n,i}^{k}$, which is a path
from $\mathbf{a}_{k}$ to $\mathbf{a}_{n}$ with $\mathbf{M}_{\alpha
}(G_{n}^{k})\leq \epsilon $. Each $d_{\alpha }-$\textbf{admissible Cauchy
sequence }corresponds to an element in $\mathcal{P}_{\alpha }(\mathbb{R}%
^{m}) $. Let
\begin{equation*}
\mathcal{D}_{\alpha }(\mathbb{R}^{m})\subset \mathcal{P}_{\alpha }(\mathbb{R}%
^{m})
\end{equation*}%
be the set of all probability measures $\mu $ which corresponds to a $%
d_{\alpha }$ admissible Cauchy sequence of probability measures. For
simplicity, we may write $\mathcal{D}_{\alpha }(\mathbb{R}^{m})$ as $%
\mathcal{D}_{\alpha }$.
\end{definition}

It is easy to see that if for each $k$, there is a partition of $\mathbf{a}%
_{k+1}=\sum_{i=1}^{N_{k}}\mathbf{a}_{k+1,i}^{(k)}$ with respect to $\mathbf{a%
}_{k}$ as sums of disjoint atomic measures and a path $G_{k+1,i}^{k}\in
Path\left( m_{i}^{(k)}\delta _{x_{i}^{(k)}},\mathbf{a}_{k+1,i}^{(k)}\right) $
for each $i=1,2,\cdots ,N_{k}$ such that
\begin{equation*}
\sum_{k=1}^{\infty }\left( \sum_{i=1}^{N_{k}}\mathbf{M}_{\alpha }\left(
G_{k+1,i}^{k}\right) \right) <+\infty ,
\end{equation*}%
then $\left\{ \mathbf{a}_{n}\right\} $ is a $d_{\alpha }$-admissible Cauchy
sequence.

Also, note that if $\mu ,\nu \in \mathcal{D}_{\alpha }(\mathbb{R}^{m})$, one
automatically has $d_{\alpha }\left( \mu ,\nu \right) <\infty $.

Before discussing properties of $\mathcal{D}_{\alpha }$, we give some
examples of elements of $\mathcal{D}_{\alpha }$ as follows. The first
example is a strengthening of corollary \ref{p=pa}.

\begin{example}
Let $\mu $ be any probability measure supported on a compact subset of $%
\mathbb{R}^{m}$. Then $\mu \in \mathcal{D}_{\alpha }$ whenever $\alpha >1-%
\frac{1}{m}$.
\end{example}

Let $B$ be a cube in $\mathbb{R}^{m}$ of side length $l$ that contains the
support of $\mu $. For each $n$, by using dyadic decomposition of $B$, we
get a family of smaller cubes $\left\{ B_{i}^{\left( n\right) }\right\}
_{i=1}^{N_{n}}$ of generation $n$ centered at $x_{i}^{\left( n\right) }$ and
of side length $\frac{l}{2^{n}}$, where $N_{n}=2^{mn}$. Then, set
\begin{equation*}
\mathbf{a}_{n}=\sum_{i=1}^{N_{n}}\mu \left( B_{i}^{\left( n\right) }\right)
\delta _{x_{i}^{(n)}.}
\end{equation*}%
For each $n$, by setting%
\begin{equation*}
\mathbf{a}_{n+1,i}^{(n)}=\sum_{x_{j}^{\left( n+1\right) }\in B_{i}^{\left(
n\right) }}\mu \left( B_{j}^{\left( n+1\right) }\right) \delta
_{x_{j}^{\left( n+1\right) }}
\end{equation*}%
for each $i=1,2,\cdots ,N_{n}$, we get a partition of $\mathbf{a}%
_{n+1}=\sum_{i=1}^{N_{n+1}}\mathbf{a}_{n+1,i}^{(n)}$ with respect to $%
\mathbf{a}_{n}$. Also, by transporting the corresponding mass $\mu \left(
B_{j}^{\left( n+1\right) }\right) $ from $x_{i}^{\left( n\right) }$ to each $%
x_{j}^{\left( n+1\right) }$, we build an obvious path $G_{n+1,i}^{n}\in
Path\left( \mu \left( B_{i}^{\left( n\right) }\right) \delta _{x_{i}^{(n)}},%
\mathbf{a}_{n+1,i}^{(n)}\right) $.

\begin{center}
\begin{tikzpicture}[style=thick]
\draw[blue,thick] (0,0) rectangle (4,4);

\foreach \x in {14,42,71,100} \foreach \y in {14,42,71,100}
{\color{red}\put(\x,\y){\circle*{4};} \color{black} \foreach \i in
{-1,1} \foreach \j in {-1,1} {\put(\x,\y){\vector(\i,\j){10};}} }
\foreach \x in {7,21,35,50,64,79,93,108} \foreach \y in
{7,21,35,50,64,79,93,108} {\color{green} \put(\x,\y){\circle*{2};} }

\end{tikzpicture}
\end{center}

Now, it is easy (see \cite[proposition 3.1]{xia1} for instance) to check
that
\begin{equation*}
\sum_{k=1}^{\infty }\left( \sum_{i=1}^{N_{k}}\mathbf{M}_{\alpha }\left(
G_{n+1,i}^{n}\right) \right) <+\infty
\end{equation*}%
whenever $\alpha >1-\frac{1}{m}$, and thus $\mu \in \mathcal{D}_{\alpha }$.

\begin{example}
\label{Cantor measure} Cantor measure
\end{example}

Let
\begin{equation*}
\mathbf{a}_{n}=\frac{1}{2^{n}}\sum_{i=1}^{2^{n}}\delta _{x_{n,i}}
\end{equation*}%
where $x_{n,i}$'s are centers of intervals of the $i^{th}$ Cantor interval
of length $\frac{1}{3^{n}}$. By transporting the mass $\frac{1}{2^{n}}$ to
each $x_{n,i}$ from the center of the interval in the previous step, we
construct a path $G_{n}^{n-1}$ from $\mathbf{a}_{n-1}$ to $\mathbf{a}_{n}$.

\begin{center}
\begin{tikzpicture}[style=thick]
\color{blue}

\foreach \x in {0,1} {\draw
(\x*4,1)--(2+\x*4,1);}
\color{red} \put(28.333,28.333){\circle*{4}} \put(141.666,28.333){\circle*{4}}

\color{blue} \foreach \x in {0,2,6,8} {\draw
(\x*2/3,1/2)--(2/3+\x*2/3,1/2);} \color{green} \foreach \x in
{9.444,47.222,122.777,160.555} {\put(\x,14.167){\circle*{2}};}

\color{green} \foreach \x in
{9.444,47.222,122.777,160.555} {\put(\x,0){\circle*{2}};}

\color{red} \put(28.333,0){\circle*{4}} \put(141.666,0){\circle*{4}}

\color{black}
\put(28.333,0){\vector(-1,0){18.888}}
\put(28.333,0){\vector(1,0){18.888}}
\put(141.666,0){\vector(-1,0){18.888}}
\put(141.666,0){\vector(1,0){18.888}}
%\caption{An example of the constructed path $G_{n}^{n-1}$ for Lebesgue measure in a box}
\end{tikzpicture}
\end{center}

%\vskip 0.3in

Now,
\begin{equation*}
\mathbf{M}_{\alpha }\left( G_{n}^{n-1}\right) =\sum_{i=1}^{2^{n}}\left(
\frac{1}{2^{n}}\right) ^{\alpha }\left( \frac{1}{3}\right) ^{n}=\left( \frac{%
2^{1-\alpha }}{3}\right) ^{n}.
\end{equation*}%
Therefore, $\left\{ \mathbf{a}_{n}\right\} $ forms a $d_{\alpha }$%
-admissible Cauchy sequence whenever $\frac{2^{1-\alpha }}{3}<1$,
that is whenever $\frac{1}{1-\alpha }>\frac{\ln 2}{\ln 3}$, which is
exactly the fractal dimension of the Cantor set. The measure
represented by this Cauchy sequence $\left\{ \mathbf{a}_{n}\right\}
$ is called the \textit{Cantor measure}. It is the usual
$\mathcal{H}^{s}\lfloor C$ where $C$ is the Cantor set and
$s=\frac{\ln 2}{\ln 3}$ is its Hausdorff dimension. This
shows that the Cantor measure is in $\mathcal{D}_{\alpha }$ whenever $\frac{1%
}{1-\alpha }>\frac{\ln 2}{\ln 3}$. Note that $\alpha $ is allowed to be
negative here.

\begin{example}
Fat Cantor measure
\end{example}

A fat Cantor set is constructed in the same way as constructing a Cantor set
except that an interval of length $\lambda $ is removed from the middle of $%
\left[ 0,1\right] $ for $\lambda \in \left( 0,1\right) $. Again, we set
\begin{equation*}
\mathbf{a}_{n}=\frac{1}{2^{n}}\sum_{i=1}^{2^{n}}\delta _{x_{n,i}}
\end{equation*}%
where $x_{n,i}$'s are centers of intervals of the $i^{th}$ Fat Cantor
interval of length $\left( \frac{1-\lambda }{2}\right) ^{n}$. By
transporting the mass $\frac{1}{2^{n}}$ to each $x_{n,i}$ from the center of
the interval in the previous step, we construct a path $G_{n}^{n-1}$ from $%
\mathbf{a}_{n-1}$ to $\mathbf{a}_{n}$. Then,
\begin{equation*}
\mathbf{M}_{\alpha }\left( G_{n}^{n-1}\right) =\sum_{i=1}^{2^{n}}\left(
\frac{1}{2^{n}}\right) ^{\alpha }\frac{1+\lambda }{4}\left( \frac{1-\lambda
}{2}\right) ^{n-1}=\frac{1+\lambda }{2\left( 1-\lambda \right) }\left(
2^{1-\alpha }p\right) ^{n},
\end{equation*}%
where $p=\frac{1-\lambda }{2}$. The sequence $\left\{ \mathbf{a}_{n}\right\}
$ forms a $d_{\alpha }$-admissible Cauchy sequence whenever $2^{1-\alpha
}p<1 $, that is whenever
\begin{equation*}
\frac{1}{1-\alpha }>-\frac{\ln 2}{\ln p}=\frac{\ln 2}{\ln 2-\ln \left(
1-\lambda \right) },
\end{equation*}%
which is the fractal dimension of the Fat Cantor set. The measure
represented by this Cauchy sequence $\left\{ \mathbf{a}_{n}\right\} $ is
called the \textit{Fat Cantor measure}. This shows that the Fat Cantor
measure is in $\mathcal{D}_{\alpha }$ whenever $\frac{1}{1-\alpha }>-\frac{%
\ln 2}{\ln p}$.

\begin{example}
Self similar measures
\end{example}

Let $A$ be any bounded self similar set in the sense that $A$ is the finite
union of sets $A_{i}$ for $i=1,\cdots ,k$ with each $A_{i}$ being a $\sigma
- $rescale of $A$. Pick any point $x^{\ast }\in A$ as the center of $A$.
Then, each copy $A_{i}^{\left( n\right) }$ of $A$ of generation $n$ has a
corresponding center $x_{i}^{\left( n\right) }$. Now, we set
\begin{equation*}
\mathbf{a}_{n}=\frac{1}{k^{n}}\sum_{i=1}^{k^{n}}\delta _{x_{i}^{\left(
n\right) }}.
\end{equation*}%
By transporting the mass $\frac{1}{k^{n}}$ to each $x_{i}^{\left( n\right) }$
from the center $x_{j}^{\left( n-1\right) }$of the set $A_{j}^{\left(
n-1\right) }$ in the previous step, we construct a path $G_{n}^{n-1}$ from $%
\mathbf{a}_{n-1}$ to $\mathbf{a}_{n}$. Note that $\left| x_{i}^{\left(
n\right) }-x_{j}^{\left( n-1\right) }\right| \lesssim \sigma ^{n-1}L$, where
$L$ is the diameter of the set $A$. Therefore,
\begin{equation*}
\mathbf{M}_{\alpha }\left( G_{n}^{n-1}\right) \approx k^{n}\left( \frac{1}{%
k^{n}}\right) ^{\alpha }\sigma ^{n-1}L=\frac{L}{\sigma }\left( k^{1-\alpha
}\sigma \right) ^{n}
\end{equation*}%
Thus,
\begin{eqnarray*}
&&\sum_{k=1}^{\infty }\mathbf{M}_{\alpha }\left( G_{n}^{n-1}\right) \text{
is finite} \\
&\Longleftrightarrow &k^{1-\alpha }\sigma <1\Longleftrightarrow \frac{1}{%
1-\alpha }>-\frac{\ln k}{\ln \sigma }.
\end{eqnarray*}%
Therefore, the sequence $\left\{ \mathbf{a}_{n}\right\} $ forms an $%
d_{\alpha }$-admissible Cauchy sequence whenever%
\begin{equation*}
\frac{1}{1-\alpha }>-\frac{\ln k}{\ln \sigma },
\end{equation*}%
which is the fractal dimension of the self similar set $A$. The measure
represented by this Cauchy sequence $\left\{ \mathbf{a}_{n}\right\} $ is
called a \textit{self similar measure}. This shows that a self similar
measure is in $\mathcal{D}_{\alpha }$ whenever $\frac{1}{1-\alpha }>-\frac{%
\ln k}{\ln \sigma }$.

\subsection{Hausdorff dimension of measures}

Let $\mathcal{H}^{s}$ denotes $s$ dimensional Hausdorff measure on $\mathbb{R%
}^{m}$ for each $s\geq 0$.

\begin{theorem}
\label{1.2} If $\mu \in \mathcal{D}_{\alpha }$ for some $\alpha <1$, then $%
\mu $ is concentrated on a subset $A$ of $\mathbb{R}^{m}$ with $\mathcal{H}^{%
\frac{1}{1-\alpha }}\left( A\right) =0$.
\end{theorem}

\begin{proof}
Since $\mu \in \mathcal{D}_{\alpha }$, it is represented by a $d_{\alpha }$%
-admissible Cauchy sequence $\left\{ \mathbf{a}_{k}\right\} $ in the form of
\begin{equation*}
\mathbf{a}_{k}=\sum_{i=1}^{N_{k}}m_{i}^{(k)}\delta _{x_{i}^{(k)}}.
\end{equation*}%
By the definition \ref{admissible} and taking a subsequence of $\left\{
\mathbf{a}_{k}\right\} $ if necessary, we have that for any $k$ and for all $%
n>k$, there exists$\text{\ a partition of }$%
\begin{equation*}
\mathbf{a}_{n}=\sum_{i=1}^{N_{k}}\mathbf{a}_{n,i}^{(k)}
\end{equation*}%
with respect to $\mathbf{a}_{k}$ and a path $G_{n,i}^{k}\in
Path(m_{i}^{(k)}\delta _{x_{i}^{(k)}},\mathbf{a}_{n,i}^{(k)})$ for each $%
i=1,2,\cdots ,N_{k}$ such that
\begin{equation*}
\sum_{i=1}^{N_{k}}d_{\alpha }(\mathbf{a}_{n,i}^{(k)},m_{i}^{(k)}\delta
_{x_{i}^{(k)}})\leq \sum_{i=1}^{N_{k}}\mathbf{M}_{\alpha }\left(
G_{n,i}^{k}\right) \leq \frac{1}{2^{k}}\text{.}
\end{equation*}%
Now, for each fixed $k$ and any $i=1,2,\cdots ,N_{k}$, the sequence $\left\{
\mathbf{a}_{n,i}^{(k)}\right\} _{n=k+1}^{\infty }$ is also a $d_{\alpha }$%
-admissible Cauchy sequence representing a positive Radon measure $\mu
_{i}^{\left( k\right) }$ of mass $m_{i}^{\left( k\right) }$. As a result,
for each fixed $k$, the measure $\mu $ can be decomposed as
\begin{equation*}
\mu =\sum_{i=1}^{N_{k}}\mu _{i}^{\left( k\right) }
\end{equation*}%
such that each $\mu _{i}^{\left( k\right) }$ has mass $m_{i}^{\left(
k\right) }$ for each $i$ and
\begin{equation*}
\sum_{i=1}^{N_{k}}d_{\alpha }(\mu _{i}^{(k)},m_{i}^{(k)}\delta
_{x_{i}^{(k)}})\leq \frac{1}{2^{k}}.
\end{equation*}

Case 1: $0<\alpha <1$. Note that from (\ref{positive_estimate_eqn}), we also
have the following estimate, for any $\nu \in \mathcal{D}_{\alpha }$ and $%
0\leq \alpha <1$,
\begin{equation*}
\nu (\mathbb{R}^{m}\setminus B(x,r))^{\alpha }r\leq d_{\alpha }(\nu ,\delta
_{x})\text{\ for any }x\in \mathbb{R}^{m}\text{\ and any }r>0.
\end{equation*}%
So if $r\geq d_{\alpha }(\nu ,\delta _{x})^{1-\alpha }$, then
\begin{equation}
\nu (\mathbb{R}^{m}\setminus B(x,r))\leq d_{\alpha }(\nu ,\delta _{x}).
\label{mu_d_alpha}
\end{equation}%
Now let $r_{i}^{(k)}=d_{\alpha }(\mu _{i}^{(k)},m_{i}^{(k)}\delta
_{x_{i}^{(k)}})^{1-\alpha }$ and $A_{k}=%
\bigcup_{i=1}^{N_{k}}B(x_{i}^{(k)},r_{i}^{\left( k\right) })$. Then
\begin{equation*}
\mathcal{H}^{\frac{1}{1-\alpha }}(A_{k})\leq \alpha \left( k\right)
\sum_{i=1}^{N_{k}}(r_{i}^{(k)})^{\frac{1}{1-\alpha }}=\alpha \left( k\right)
\sum_{i=1}^{N_{k}}d_{\alpha }(\mu _{i}^{(k)},m_{i}^{(k)}\delta
_{x_{i}^{(k)}})\leq \frac{\alpha \left( k\right) }{2^{k}},
\end{equation*}%
where $\alpha \left( k\right) $ is the constant
\begin{equation}
\alpha \left( k\right) =\frac{\pi ^{\frac{k}{2}}}{\Gamma \left( \frac{k}{2}%
+1\right) }.  \label{Hausdorff constant}
\end{equation}%
Moreover,
\begin{eqnarray*}
\mu (\mathbb{R}^{m}\setminus A_{k}) &\leq &\sum_{i=1}^{N_{k}}\mu
_{i}^{\left( k\right) }(\mathbb{R}^{m}\setminus B(x_{i}^{(k)},r_{i}^{\left(
k\right) })) \\
&\overset{\text{ by }(\ref{mu_d_alpha})}{\leq }&\sum_{i=1}^{N_{k}}d_{\alpha
}(\mu _{i}^{(k)},m_{i}^{(k)}\delta _{x_{i}^{(k)}})\leq \frac{1}{2^{k}}.
\end{eqnarray*}

We set $A:=\bigcup_{h}\left( \bigcap_{k>h}A_{k}\right) $. Then, for each $h$%
,
\begin{eqnarray*}
\mu \left( \mathbb{R}^{m}\setminus A\right) &\leq &\mu (\mathbb{R}%
^{m}\setminus \left( \bigcap_{k>h}A_{k}\right) )=\mu (\bigcup_{k>h}(\mathbb{R%
}^{m}\setminus A_{k})) \\
&\leq &\sum_{k>h}\mu (\mathbb{R}^{m}\setminus A_{k})\leq \sum_{k>h}\frac{1}{%
2^{k}}=\frac{1}{2^{h}}.
\end{eqnarray*}%
This implies that $\mu \left( \mathbb{R}^{m}\setminus A\right) =0$. On the
other hand, since
\begin{equation*}
\mathcal{H}^{\frac{1}{1-\alpha }}(\bigcap_{k>h}A_{k})\leq \mathcal{H}^{\frac{%
1}{1-\alpha }}(A_{k})\leq \frac{\alpha \left( k\right) }{2^{k}}\rightarrow 0%
\text{,}
\end{equation*}%
we have $\mathcal{H}^{\frac{1}{1-\alpha }}(\bigcap_{k>h}A_{k})=0$ for all $%
h\in \mathbb{N}$ and thus $\mathcal{H}^{\frac{1}{1-\alpha }}(A)=0$. Thus $%
\mu $ is concentrated on a $\frac{1}{1-\alpha }$-negligible set $A$.

Case 2: $\alpha \leq 0$. We shall denote
\begin{equation*}
r_{i}^{(k)}=\frac{d_{\alpha }(m_{i}^{(k)}\delta _{x_{i}^{(k)}},\mu
_{i}^{(k)})}{(m_{i}^{(k)})^{\alpha }}
\end{equation*}%
and $A_{k}=\bigcup_{i=1}^{N_{k}}B(x_{i}^{(k)},r_{i}^{\left( k\right) }).$

Suppose the atomic measure $\mathbf{a}_{n,i}^{(k)}$ is expressed as
\begin{equation*}
\mathbf{a}_{n,i}^{(k)}=\sum_{j=p_{i}}^{q_{i}}m_{j}^{(n)}\delta
_{x_{j}^{(n)}}.
\end{equation*}%
Then, by corollary \ref{negative_estimate}, for every $x_{j}^{(n)}$, $%
j=p_{i},...,q_{i}$ we have
\begin{equation*}
d_{\alpha }(m_{i}^{(k)}\delta _{x_{i}^{(k)}},\mathbf{a}_{n,i}^{(k)})\geq
(m_{i}^{(k)})^{\alpha }\left| x_{i}^{(k)}-x_{j}^{(n)}\right| .
\end{equation*}%
Therefore $\left| x_{i}^{(k)}-x_{j}^{(n)}\right| \leq r_{i}^{(k)}$, so $%
x_{j}^{(n)}\in B(x_{i}^{(k)},r_{i}^{\left( k\right) }).$ Thus $\mu (\mathbb{R%
}^{m}\setminus A_{k})=0.$

Moreover
\begin{eqnarray*}
\mathcal{H}^{\frac{1}{1-\alpha }}(A_{k}) &\leq
&\sum_{i=1}^{N_{k}}(r_{i}^{(k)})^{\frac{1}{1-\alpha }}=%
\sum_{i=1}^{N_{k}}(d_{\alpha }(m_{i}^{(k)}\delta _{x_{i}^{(k)}},\mu
_{i}^{(k)}))^{\frac{1}{1-\alpha }}(m_{i}^{(k)})^{\frac{-\alpha }{1-\alpha }}
\\
&\leq &(\sum_{i=1}^{N_{k}}(d_{\alpha }(m_{i}^{(k)}\delta _{x_{i}^{(k)}},\mu
_{i}^{(k)}))^{\frac{1}{1-\alpha }\cdot (1-\alpha )})^{\frac{1}{1-\alpha }%
}(\sum_{i=1}^{N_{k}}(m_{i}^{(k)})^{\frac{-\alpha }{1-\alpha }\cdot \frac{%
1-\alpha }{-\alpha }})^{\frac{-\alpha }{1-\alpha }} \\
&=&d_{\alpha }(\mu ,\mathbf{a}_{k})^{\frac{1}{1-\alpha }}.
\end{eqnarray*}

For $A=\bigcap A_{k}$. Then
\begin{equation*}
\mu (\mathbb{R}^{m}\setminus A)=\mu (\bigcup_{k}(\mathbb{R}^{m}\setminus
A_{k}))\leq \sum_{k}\mu (\mathbb{R}^{m}\setminus A_{k})=0.
\end{equation*}%
Moreover, we have $\mathcal{H}^{\frac{1}{1-\alpha }}(A)\leq \mathcal{H}^{%
\frac{1}{1-\alpha }}(A_{k})\leq d_{\alpha }(\mu ,\mathbf{a}_{k})^{\frac{1}{%
1-\alpha }}\rightarrow 0$ as $k\rightarrow \infty $. Thus $\mu $ is
concentrated on a $\frac{1}{1-\alpha }$-negligible set $A$.
\end{proof}

\begin{definition}
For any probability measure $\mu $ on $\mathbb{R}^{m}$, the Hausdorff
dimension of $\mu $ is defined to be
\begin{equation*}
\dim _{H}\left( \mu \right) =\inf \left\{ \dim _{H}\left( A\right) :\mu
\left( \mathbb{R}^{m}\backslash A\right) =0\right\} ,
\end{equation*}%
where $\dim _{H}(A)$ is the Hausdorff dimension of a set $A$.
\end{definition}

Thus, by Theorem \ref{1.2}, we have

\begin{corollary}
For any $\alpha <1$ and any $\mu \in \mathcal{D}_{\alpha }(\mathbb{R}^{m})$,
we have
\begin{equation*}
\dim _{H}(\mu )\leq \frac{1}{1-\alpha }.
\end{equation*}
\end{corollary}

\begin{lemma}
\label{cantor}Let $\mu $ be the Cantor measure as defined in example \ref%
{cantor}. Then, $\dim _{H}\left( \mu \right) =\frac{\ln 2}{\ln 3}$
\end{lemma}

\begin{proof}
Since $\mu $ is clearly concentrated on the Cantor set whose Hausdorff
dimension is $\frac{\ln 2}{\ln 3}$, we have $\dim _{H}\left( \mu \right)
\leq \frac{\ln 2}{\ln 3}$. Therefore, to show $\dim _{H}\left( \mu \right) =%
\frac{\ln 2}{\ln 3}$, all we need to show is
\begin{equation*}
\dim _{H}\left( \mu \right) \geq \frac{\ln 2}{\ln 3}:=s.
\end{equation*}%
For this, it is sufficient to show that
\begin{equation*}
H^{s}\left( A\right) \geq 1
\end{equation*}%
whenever $\mu $ is concentrated on a set $A$. Let $\left\{ C_{i}\right\} $
be any collection of sets that covers $A$. We want to show that
\begin{equation}
\sum_{i}\alpha \left( s\right) \left( \frac{diam(C_{i})}{2}\right)
^{s}\geq 1, \label{Hausdorff_Cantor}
\end{equation}%
where $\alpha \left( s\right) $ is the constant as given in (\ref{Hausdorff
constant}). Without losing generality, we may assume that $C_{i}$ are closed
intervals. Also, if $C_{i}\cap C_{j}\neq \emptyset $ for some $i\neq j$,
then
\begin{equation*}
\left( diam\left( C_{i}\cup C_{j}\right) \right) ^{s}\leq \left(
diam(C_{i})+diam(C_{j})\right) ^{s}\leq \left( diam(C_{i})\right)
^{s}+\left( diam(C_{j})\right) ^{s}
\end{equation*}%
as the function $x^{s}$ is concave on $\left[ 0,1\right] $. Therefore, to
prove the inequality (\ref{Hausdorff_Cantor}), we may replace $C_{i}$ and $%
C_{j}$ by $C_{i}\cup C_{j}$. Thus, without losing generality, we may assume
that $\left\{ C_{i}\right\} $ are disjoint closed intervals $C_{i}=\left[
a_{i},b_{i}\right] $ with the order
\begin{equation*}
0\leq a_{1}<b_{1}<a_{2}<b_{2}<\cdots <a_{n}<b_{n}<\cdots \leq 1.
\end{equation*}%
Note that for each $i,$%
\begin{equation*}
\mu \left( \left( b_{i},a_{i+1}\right) \right) =0
\end{equation*}%
as $\left( b_{i},a_{i+1}\right) \cap A=\emptyset $ and $\mu $ is
concentrated on $A$. Therefore, $\mathbf{a}_{n}\left( \left(
b_{i},a_{i+1}\right) \right) =0$ when $n$ is large enough. This implies that
the Cantor set is disjoint with $\left( b_{i},a_{i+1}\right) $. Thus, the
Cantor set is also covered by $\left\{ C_{i}\right\} $. As a result, $%
H^{s}\left( \cup C_{i}\right) \geq H^{s}\left( \text{Cantor set}\right) =1$.
This shows the inequality (\ref{Hausdorff_Cantor}) and hence $\dim
_{H}\left( \mu \right) =s$ as desired.
\end{proof}

\subsection{Minkowski dimension of measures}

A \textit{nested collection}
\begin{equation*}
\mathcal{F}=\left\{ Q_{i}^{n}:i=1,2,\cdots ,N_{n}\text{ and }n=1,2,\cdots
\right\}
\end{equation*}
is a collection of Borel subsets of $\mathbb{R}^{m}$ with the following
properties:

\begin{enumerate}
\item for each $Q_{i}^{n}$, its diameter
\begin{equation}
C_{1}\sigma ^{n}\leq diam\left( Q_{i}^{n}\right) \leq C_{2}\sigma ^{n}
\label{condition_1}
\end{equation}%
for some constants $C_{2}\geq C_{1}>0$ and some $\sigma \in \left(
0,1\right) $.

\item for any $k,l,i,j$ with $l\leq k$, either $Q_{i}^{k}\cap
Q_{j}^{l}=\emptyset $ or $Q_{i}^{k}\subseteq Q_{j}^{l};$

\item for each $Q_{j}^{n+1}$ there exists exactly one $Q_{i}^{n}$ (parent of
$Q_{j}^{n+1}$) such that $Q_{j}^{n+1}\subseteq Q_{i}^{n}$;

\item for each $Q_{i}^{n}$ there exists at least one $Q_{j}^{n+1}$ (child of
$Q_{i}^{n}$) such that $Q_{j}^{n+1}\subseteq Q_{i}^{n}$;
\end{enumerate}

Each $Q_{i}^{n}$ is called a cube of generation $n$ in $\mathcal{F}$. If two
different cubes $Q_{i}^{n}$ and $Q_{j}^{n}$ of generation $n$ have the same
parent, then they are called \textit{brothers} to each other.

A typical example of a nested collection includes collections of standard
cubes in $\mathbb{R}^{m}$. That is, let $Q$ be a cube in $\mathbb{R}^{m}$ of
side length $L$ and $k$ is a fixed natural number. Then, we evenly split $Q$
into $k^{m}$ cubes of side length $\frac{L}{k}$. Pick some (or all) of these
cubes of generation $1$ to form a collection $\mathcal{Q}_{1}$. We may then
evenly split each of the cubes of generation $1$ in $\mathcal{Q}_{1}$ into
cubes of side length $\frac{L}{k^{2}}$ to get cubes of generation $2$. Pick
at least one cube of generation $2$ from each cube in $\mathcal{Q}_{1}$, we
get a collection $\mathcal{Q}_{2}$ consisting of cubes of generation $2$.
Then, we may continue this process for each $n$ to get a collection $%
\mathcal{Q}_{n}$ consisting some cubes of side length $\frac{L}{k^{n}}$. The
union of all $\mathcal{Q}_{n}$ is clearly a nested collection.

\begin{definition}
For any nested collection $\mathcal{F}$, we define its Minkowski dimension
\begin{equation}
\dim _{M}\left( \mathcal{F}\right) :=\lim_{n\rightarrow \infty }\frac{\log
\left( N_{n}\right) }{\log \left( \frac{1}{\sigma ^{n}}\right) }
\label{cubical_dim}
\end{equation}%
provided the limit exists, where $N_{n}$ is the total number of cubes of
generation $n$.
\end{definition}

By (\ref{condition_1}), we have
\begin{equation*}
\dim _{M}\left( \mathcal{F}\right) \leq m
\end{equation*}%
for any nested collection $\mathcal{F}$ consisting Borel subsets of $\mathbb{%
R}^{m}$.

\begin{definition}
A Radon measure $\mu $ is said to be \textit{concentrated} on a nested
collection $\mathcal{F}$ if for each $n$,
\begin{equation*}
\mu \left( \mathbb{R}^{m}\setminus \left(
\bigcup_{i=1}^{N_{n}}Q_{i}^{n}\right) \right) =0.
\end{equation*}
\end{definition}

Now, we inductively define the centers for cubes in $\mathcal{F}$ as
follows. For any cube $Q_{i}^{n}$ in $\mathcal{F}$, if it has more than one
child, then one may pick any point $p$ in $Q_{i}^{n}$, and call it the
center of $Q_{i}^{n}$. If $Q_{i}^{n}$ has only one child, then we pick the
center of the child as the center of $Q_{i}^{n}$. In case that we have an
infinite sequence of cubes $\left\{ Q_{i}^{n}\right\} $ such that each cube
is the only child of the previous one, then the intersection of the closures
of these cubes contains at least one point, and we will call this point the
center of each cube in the sequence. It is easy to see that we have defined
a center for each cube in $\mathcal{F}$ by this process.

Now, for each $Q_{i}^{n}$ in $\mathcal{F}$, let
\begin{equation*}
l\left( Q_{i}^{n}\right)
\end{equation*}%
be the distance from the center of $Q_{i}^{n}$ to the center of its parent $%
Q_{j}^{n-1}$. Then, by definition of centers and (\ref{condition_1}), we
have
\begin{equation}
l\left( Q_{i}^{n}\right) =\left\{
\begin{array}{cc}
0, & \text{if }Q_{i}^{n}\text{ has no brothers,} \\
\leq C_{2}\sigma ^{n-1} & \text{otherwise.}%
\end{array}%
\right.  \label{center_estimate}
\end{equation}%
We now state our key lemma as follows:

\begin{lemma}
\label{key_lemma}Suppose $\mu $ is a probability measure concentrated on a
nested collection $\mathcal{F}$ with
\begin{equation*}
\dim _{M}\left( \mathcal{F}\right) <\frac{1}{1-\alpha }
\end{equation*}%
for some $\alpha <1$. If for each $n,$
\begin{equation}
\sum_{i=1}^{N_{n}}\mu \left( Q_{i}^{n}\right) ^{\alpha }l\left(
Q_{i}^{n}\right) \leq C\left( N_{n}\right) ^{1-\alpha }\sigma ^{n}
\label{key_condition}
\end{equation}%
for some constant $C$, then $\mu \in \mathcal{D}_{\alpha }$.
\end{lemma}

\begin{proof}
Let
\begin{equation*}
\beta :=\dim _{M}\left( \mathcal{F}\right) <\frac{1}{1-\alpha }.
\end{equation*}%
Then, by (\ref{cubical_dim}), when $n$ is large enough, we have
\begin{equation*}
N_{n}\leq \sigma ^{-n\beta }.
\end{equation*}%
Now, we construct a $d_{\alpha }-$ admissible Cauchy sequence $\left\{
\mathbf{a}_{n}\right\} $ for $\mu $ as follows. For each $n$, let
\begin{equation*}
\mathbf{a}_{n}=\sum_{i=1}^{N_{n}}\mu (Q_{i}^{n})\delta _{x_{i}^{(n)}},
\end{equation*}%
where $x_{i}^{\left( n\right) }$ is the center of $Q_{i}^{n}$ and $%
m_{i}^{\left( n\right) }=\mu (Q_{i}^{n})$ denotes the total mass of $\mu $
on $Q_{i}^{n}$. Also, we can construct a partition of $\mathbf{a}_{n}$ with
respect to $\mathbf{a}_{n-1}$ by grouping together masses that are located
in centers of children of the same cube. That is, for each $j=1,2,\cdots
,N_{n-1}$, let
\begin{equation*}
\mathbf{a}_{n-1,j}^{\left( n\right) }=\sum \mu (Q_{i}^{n})\delta
_{x_{i}^{(n)}}
\end{equation*}%
where the summation is over all $i^{\prime }$s such that $Q_{i}^{n}$ is a
child of $Q_{j}^{n-1}$. For each $Q_{i}^{n}$, we may transport the mass $%
m_{i}^{\left( n\right) }$ from the center of its parent to its center, and
thus construct a path $G_{n}^{n-1}$ from $\mathbf{a}_{n-1}$ to $\mathbf{a}%
_{n}$. Then, by (\ref{key_condition}),
\begin{eqnarray*}
\mathbf{M}_{\alpha }\left( G_{n}^{n-1}\right) &\leq
&\sum_{i=1}^{N_{n}}\left( \mu (Q_{i}^{n})\right) ^{\alpha }l\left(
Q_{i}^{n}\right) \\
&\leq &C\left( N_{n}\right) ^{1-\alpha }\sigma ^{n} \\
&\leq &C\left( \sigma ^{-n\beta }\right) ^{1-\alpha }\sigma ^{n}\text{ when }%
n\text{ is large enough} \\
&=&Cb^{n}
\end{eqnarray*}%
for $b=\sigma ^{1-\beta \left( 1-\alpha \right) }\in \left( 0,1\right) $. As
a result, $\left\{ \mathbf{a}_{n}\right\} $ is a $d_{\alpha }$-admissible
Cauchy sequence representing $\mu $ and thus $\mu \in \mathcal{D}_{\alpha }$.
\end{proof}

\begin{definition}
For any Radon measure $\mu $, we define the Minkowski dimension of the
measure $\mu $ to be
\begin{equation*}
\dim _{M}\left( \mu \right) :=\inf \left\{ \dim _{M}\left( \mathcal{F}%
\right) \right\}
\end{equation*}%
where the infimum is over all nested collection $\mathcal{F}$ that $\mu $ is
concentrated on.
\end{definition}

\begin{theorem}
\label{3.1} Suppose $\mu $ is a probability measure with $\dim _{M}(\mu )<%
\frac{1}{1-\alpha }$ for some $0\leq \alpha <1$, then $\mu \in \mathcal{D}%
_{\alpha }$.
\end{theorem}

\begin{proof}
Since $\dim _{M}(\mu )<\frac{1}{1-\alpha }$, $\mu $ is concentrated on a
nested collection $\mathcal{F}$ with
\begin{equation*}
\dim _{M}\left( \mathcal{F}\right) <\frac{1}{1-\alpha }.
\end{equation*}%
When $\alpha \geq 0$, by \textit{H\"{o}lder inequality }and (\ref%
{center_estimate}) we have\textit{\ }%
\begin{eqnarray*}
&&\sum_{i=1}^{N_{n}}\mu \left( Q_{i}^{n}\right) ^{\alpha }l\left(
Q_{i}^{n}\right) \\
&\leq &\left( \sum_{i=1}^{N_{n}}\mu \left( Q_{i}^{n}\right) \right) ^{\alpha
}(\sum_{i=1}^{N_{n}}1)^{1-\alpha }C_{2}\sigma ^{n-1} \\
&=&C_{2}N_{n}^{1-\alpha }\sigma ^{n-1}\text{, \textit{because }}%
\sum_{i=1}^{N_{n}}\mu \left( Q_{i}^{n}\right) =1
\end{eqnarray*}%
Thus, by lemma \ref{key_lemma}, we have $\mu \in \mathcal{D}_{\alpha }$.
\end{proof}

\subsection{Evenly concentrated measures}

Now, we aim at achieving a similar result as in theorem \ref{3.1} for the
case $\alpha <0$. To do it, we introduce the following definition:

\begin{definition}
A Radon measure $\mu $ is \textit{evenly concentrated} on a nested
collection $\mathcal{F}$ if for each cube $Q_{i}^{n}$ of generation $n$ in $%
\mathcal{F}$, either $Q_{i}^{n}$ has no brothers or $\mu \left(
Q_{i}^{n}\right) \geq \frac{\lambda }{N_{n}}$ for some constant $\lambda >0$.
\end{definition}

Here, $Q_{i}^{n}$ has no brothers means that the parent of $Q_{i}^{n}$ has
only one child, namely $Q_{i}^{n}$ itself.

Now, we provide some examples of a measure $\mu $ that is evenly
concentrated on some nested collection.
\begin{figure}[h!]
\centering \includegraphics[width=\textwidth]{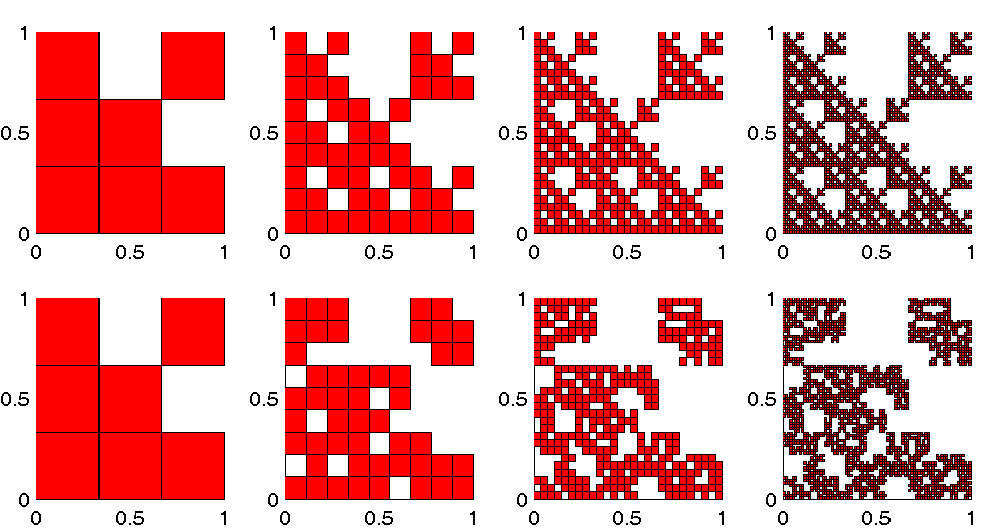}
\caption{Here, $k=3, m=2$ and $h=7$. Boxes in the first row is chosen
determinately, and boxes in the second row is picked randomly.}
\end{figure}

\begin{example}
Suppose $k$ is a fixed natural number and $h$ is a natural number no more
than $k^{m}$. Let $B_{0}$ be a fixed box in $\mathbb{R}^{m}$ and we set $%
\mathcal{B}_{1}=\left\{ B_{0}\right\} $. Now, for each $n\geq 1$, and for
each $B\in \mathcal{B}_{n}$, we split $B$ evenly into $k^{m}$ many smaller
boxes of length size $\sigma ^{n}$ and then pick $h$ smaller boxes from them
either determinately or randomly. Let $\mathcal{B}_{n+1}$ be the collections
of all smaller boxes as picked above. Thus, for each $n$, the cardinality of
$\mathcal{B}_{n+1}$ is $N_{n+1}=hN_{n}=h^{n}$ for each $n$. Let
\begin{equation*}
\mathbf{a}_{n}=\sum_{B\in \mathcal{B}_{n}}\frac{1}{N_{n}}\delta _{x_{B}}
\end{equation*}%
where $x_{B}$ is the center of the box $B\in \mathcal{B}_{n}$. Then, $%
\left\{ \mathbf{a}_{n}\right\} $ is a sequence of probability measures. For
each $B\in \mathcal{B}_{n}$, we transport the mass $\frac{1}{N_{n}}$ from
the center $x_{B^{\ast }}$ of $B$'s parent $B^{\ast }\in \mathcal{B}_{n-1}$
to $x_{B}$. In this way, we construct a path $G_{n}^{n-1}$ from $\mathbf{a}%
_{n-1}$ to $\mathbf{a}_{n}$. For each $\alpha $,
\begin{equation*}
\mathbf{M}_{\alpha }\left( G_{n}^{n-1}\right) \simeq
\sum_{i=1}^{N_{n}}\left( \frac{1}{N_{n}}\right) ^{\alpha }\left( \frac{1}{k}%
\right) ^{n}=\frac{\left( N_{n}\right) ^{1-\alpha }}{k^{n}}=\frac{1}{%
h^{1-\alpha }}\left( \frac{h^{1-\alpha }}{k}\right) ^{n}.
\end{equation*}%
Thus, $\left\{ \mathbf{a}_{n}\right\} \in \mathcal{D}_{\alpha }$ if $\frac{%
h^{1-\alpha }}{k}<1$ i.e. $\frac{1}{1-\alpha }>\log _{k}^{h}$, which is the
Minkowski dimension of the nested collection $\mathcal{B}=\bigcup_{n}%
\mathcal{B}_{n}$. Hence, $\left\{ \mathbf{a}_{n}\right\} $ is a $d_{\alpha }$%
-Cauchy sequence representing a probability measure $\mu $. Moreover, for
each $j>0,$
\begin{equation*}
\mathbf{a}_{n+j}\left( B\right) =\sum_{\substack{ B^{\ast }\in \mathcal{B}%
_{n+j}  \\ B^{\ast }\subset B}}\frac{1}{N_{n+j}}\delta _{x_{B}^{\ast }}=%
\frac{h^{j}}{N_{n+j}}=\frac{1}{N_{n}}=\mathbf{a}_{n}\left( B\right) .
\end{equation*}

Hence, $\mu $ is evenly concentrated on the nested collection $\mathcal{B}$
with the property that $\mu \left( B\right) =\frac{1}{N_{n}}$ for each $B\in
\mathcal{B}_{n}.$
\end{example}

\begin{example}
Now, we modify the previous example as follows.
\begin{figure}[h!]
\centering \includegraphics[width=\textwidth]{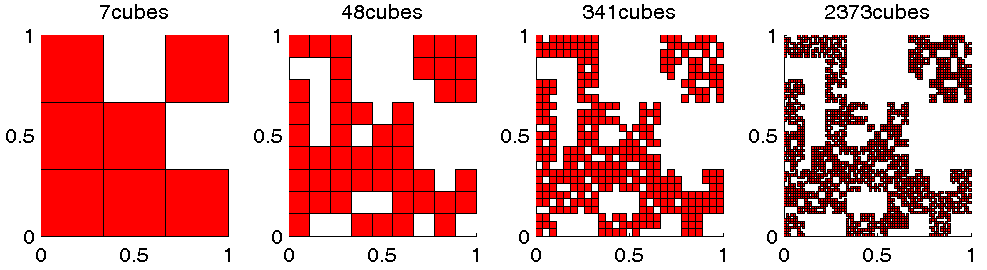}
\caption{Here, $k=3,m=2$ and $h=7$. We divide every cube into 9 smaller
cubes and then pick about 7 ones from them with some randomness. For any
particular cube, if the number of cubes picked at the current stage is less
than (or respectively greater than) the desired one (i.e. 7), more (or
respectively less) smaller cubes inside it will be picked for the next
generation.}
\end{figure}
For any $B\in \mathcal{B}_{n}$, suppose
\begin{equation*}
\lambda _{2}\leq \frac{card\left\{ B^{\ast }\in \mathcal{B}_{n+j}:B^{\ast
}\subset B\right\} }{h^{j}}\leq \lambda _{1}
\end{equation*}%
whenever $j$ is large enough for some constants $0<\lambda _{2}\leq \lambda
_{1}$. Define $\left\{ \mathbf{a}_{n}\right\} ,$ $\mathcal{B}$ and $\mu $ as
before. Thus, when $j$ is large enough
\begin{eqnarray*}
N_{n+j} &=&\sum_{B\in \mathcal{B}_{n}}card\left\{ B^{\ast }\in \mathcal{B}%
_{n+j}:B^{\ast }\subset B\right\} \\
&\leq &\sum_{B\in \mathcal{B}_{n}}\lambda _{1}h^{j}=\lambda _{1}N_{n}h^{j}.
\end{eqnarray*}

Then,
\begin{eqnarray*}
\mu \left( B\right) &=&\lim_{j}\mathbf{a}_{n+j}\left( B\right) \\
&=&\lim_{j}\frac{card\left\{ B^{\ast }\in \mathcal{B}_{n+j}:B^{\ast }\subset
B\right\} }{N_{n+j}} \\
&\geq &\lim_{j}\frac{\lambda _{2}h^{j}}{\lambda _{1}N_{n}h^{j}}=\frac{%
\lambda }{N_{n}}\text{ for }\lambda =\frac{\lambda _{2}}{\lambda _{1}}\text{.%
}
\end{eqnarray*}%
Therefore, in this case, $\mu $ is still evenly concentrated on a nested
collection $\mathcal{B}$.
\end{example}

Recall that a nonnegative Borel measure $\mu $ is Ahlfors regular of
dimension $d$ if there exists a constant $C>0$ such that
\begin{equation*}
C^{-1}r^{d}\leq \mu \left( B\left( x,r\right) \right) \leq Cr^{d}
\end{equation*}%
whenever $0<r\leq diam\left( spt\left( \mu \right) \right) $ and $x$ lies in
the support $spt\left( \mu \right) $ of $\mu $.

\begin{example}
If $\mu $ is an Ahlfors regular measure whose support is a box $\tilde{B}$
in $\mathbb{R}^{m}$ , then $\mu $ is evenly concentrated on a nested
collection.

\begin{proof}
By splitting the box $\tilde{B}$ into smaller boxes of length size $\sigma
^{n}$, we get a collection $\mathcal{B}_{n}$ of disjoint boxes. For any $%
B\in \mathcal{B}_{n}$, we have
\begin{equation*}
a\left( \sigma ^{n}\right) ^{d}\leq \mu \left( B\right) \leq A\left( \sigma
^{n}\right) ^{d}
\end{equation*}%
for some $d$ and for some suitable constants $a,A>0$. Thus,
\begin{equation*}
\mu \left( \tilde{B}\right) =\sum_{B\in \mathcal{B}_{n}}\mu \left( B\right)
\leq \sum_{B\in \mathcal{B}_{n}}A\left( \sigma ^{n}\right) ^{d}\leq
AN_{n}\sigma ^{nd}
\end{equation*}%
where $N_{n}$ is the cardinality of the collection $\mathcal{B}_{n}$.
Therefore, for every $B\in \mathcal{B}_{n}$
\begin{equation*}
\mu \left( B\right) \geq a\sigma ^{nd}\geq a\frac{\mu \left( \tilde{B}%
\right) }{AN_{n}}=\frac{\lambda }{N_{n}}
\end{equation*}%
where $\lambda =\frac{a}{A}\mu \left( \tilde{B}\right) >0$. Hence, $\mu $ is
evenly concentrated on a nested collection.
\end{proof}
\end{example}

\begin{definition}
For any Radon measure $\mu $, we define
\begin{equation*}
\dim _{U}\left( \mu \right) :=\inf \left\{ \dim _{M}\left( \mathcal{F}%
\right) \right\}
\end{equation*}%
where the infimum is over all nested collection $\mathcal{F}$ that $\mu $ is
evenly concentrated on.
\end{definition}

Obviously,
\begin{equation*}
\dim _{M}\left( \mu \right) \leq \dim _{U}\left( \mu \right) .
\end{equation*}
Now, we have the following theorem for any $\alpha<1$.

\begin{theorem}
\label{negative}Suppose $\mu $ is a probability measure with $\dim _{U}(\mu
)<\frac{1}{1-\alpha }$ for some $\alpha <1$, then $\mu \in \mathcal{D}%
_{\alpha }$.
\end{theorem}

\begin{proof}
Since $\dim _{M}\left( \mu \right) \leq \dim _{U}\left( \mu \right) $, by
theorem \ref{3.1}, $\mu \in \mathcal{D}_{\alpha }$ whenever $0\leq \alpha <1$%
. Thus, we only need to consider the case that $\alpha <0$. Since $\dim
_{U}(\mu )<\frac{1}{1-\alpha }$, $\mu $ is evenly concentrated on a nested
collection $\mathcal{F}$ with
\begin{equation*}
\dim _{M}\left( \mathcal{F}\right) <\frac{1}{1-\alpha }.
\end{equation*}%
For each cube $Q_{i}^{n}$ in $\mathcal{F}$, if $Q_{i}^{n}$ has no brothers,
then $l\left( Q_{i}^{n}\right) =0$. If $Q_{i}^{n}$ has more than one
brother, then
\begin{equation*}
\mu \left( Q_{i}^{n}\right) \geq \frac{\lambda }{N_{n}}\text{ and }l\left(
Q_{i}^{n}\right) \leq C_{2}\sigma ^{n-1}\text{.}
\end{equation*}%
Thus, when $\alpha \leq 0$, we always have
\begin{equation*}
\mu \left( Q_{i}^{n}\right) ^{\alpha }l\left( Q_{i}^{n}\right) \leq \left(
\frac{\lambda }{N_{n}}\right) ^{\alpha }C_{2}\sigma ^{n-1}.
\end{equation*}%
Therefore, \
\begin{equation*}
\sum_{i=1}^{N_{n}}\mu \left( Q_{i}^{n}\right) ^{\alpha }l\left(
Q_{i}^{n}\right) \leq \frac{C_{2}\lambda ^{\alpha }}{\sigma }\left(
N_{n}\right) ^{1-\alpha }\sigma ^{n}.
\end{equation*}%
By lemma \ref{key_lemma}, we still have $\mu \in \mathcal{D}_{\alpha }$ when
$\alpha \leq 0$.
\end{proof}

\subsection{Comparison of dimensions of measures}

\begin{definition}
For any probability measure $\mu $, we define the transport dimension of $%
\mu $ to be
\begin{equation*}
\dim _{T}\left( \mu \right) =\inf_{\alpha <1}\left\{ \frac{1}{1-\alpha }:\mu
\in \mathcal{D}_{\alpha }(\mathbb{R}^{m})\right\} .
\end{equation*}
\end{definition}

\begin{theorem}
\label{main theorem} Let $\mu $ be any positive probability measure on $%
\mathbb{R}^{m}$. Then, the following bounds hold
\begin{equation*}
\dim _{H}(\mu )\leq \dim _{T}(\mu )\leq max\{\dim _{M}(\mu ),1\}.
\end{equation*}%
Moreover, we also have
\begin{equation*}
\dim _{H}(\mu )\leq \dim _{T}(\mu )\leq \dim _{U}(\mu ).
\end{equation*}
\end{theorem}

\begin{proof}
Assume that $\dim _{T}(\mu )>\max \{\dim _{M}(\mu ),1\}$ (or $\dim _{U}(\mu
) $, respectively ). Then, we may pick an $\alpha <1$ so that
\begin{equation*}
\dim _{T}(\mu )>\frac{1}{1-\alpha }>\max \{\dim _{M}(\mu ),1\}\text{ }\left(
\text{or }\dim _{U}\left( \mu \right) \text{ respectively}\right) .
\end{equation*}%
By Theorem \ref{3.1} (or Theorem \ref{negative}), $\mu \in \mathcal{D}%
_{\alpha }(\mathbb{R}^{m})$. This is a contradiction to the definition of $%
\dim_{T} \left( \mu \right) $.

On the other hand, by Theorem \ref{1.2}, $\dim _{H}(\mu )\leq \frac{1}{%
1-\alpha }$ whenever $\mu \in \mathcal{D}_{\alpha }(\mathbb{R}^{m})$. Thus, $%
\dim _{H}(\mu )\leq \dim _{T}(\mu )$.
\end{proof}

\begin{example}
\thinspace Let $\mu $ be the Cantor measure as in Example \ref{Cantor
measure}. Let $Q_{i}^{n}$ denotes the $i^{th}$ Cantor interval of length $%
\frac{1}{^{3^{n}}}$. Then, clearly $\mu $ is evenly concentrated on the
nested collection $\mathcal{C}=\left\{ Q_{i}^{n}:i=,2,\cdots ,2^{n}\right\}
_{n=1}^{\infty }$. Note that
\begin{equation*}
\dim _{M}\left( \mathcal{C}\right) =\lim_{n\rightarrow \infty }\frac{\log
\left( N_{n}\right) }{\log \left( 3^{n}\right) }=\lim_{n\rightarrow \infty }%
\frac{\log \left( 2^{n}\right) }{\log \left( 3^{n}\right) }=\frac{\ln 2}{\ln
3}.
\end{equation*}%
On the other hand, by lemma \ref{cantor}, $\dim _{H}\left( \mu \right) =%
\frac{\ln 2}{\ln 3}$. Therefore, by theorem \ref{main theorem},
\begin{equation*}
\dim _{T}\left( \mu \right) =\frac{\ln 2}{\ln 3}.
\end{equation*}
\end{example}

\section{The Dimensional Distance between probability measures}

For any $\alpha <1$, let
\begin{equation*}
\mathcal{S}_{\alpha }\left( \mathbb{R}^{m}\right) =\left\{ \Lambda \left(
\mu -\nu \right) :\Lambda \geq 0,\mu ,\nu \in \mathcal{D}_{\alpha }(\mathbb{R%
}^{m})\right\}
\end{equation*}%
be a collection of signed measures. Clearly, $\mathcal{S}_{\alpha
_{1}}\left( \mathbb{R}^{m}\right) \subseteq \mathcal{S}_{\alpha _{2}}\left(
\mathbb{R}^{m}\right) $ if $\alpha _{1}\leq \alpha _{2}$.

\begin{definition}
\label{Def_D} For any two probability measures $\mu ,\nu $ on $\mathbb{R}%
^{m} $, we define%
\begin{equation*}
D(\mu ,\nu ):=\inf_{\alpha <1}\{\frac{1}{1-\alpha }:\mu -\nu \in \mathcal{S}%
_{\alpha }\left( \mathbb{R}^{m}\right) \}.
\end{equation*}
\end{definition}

\begin{proposition}
\label{pseudometric} $D\left( \mu ,\nu \right) \,$ is a pseudometric%
\footnote[1]{%
A pseudometric $D$ means that it is nonnegative, symmetric, satisfies the
triangle inequality, and $D\left( \mu ,\mu \right) =0$. But $D\left( \mu
,\nu \right) =0$ does not imply $\mu =\nu $.} on the space of probability
measures on $\mathbb{R}^{m}$.
\end{proposition}

\begin{proof}
It is clear that $D$ is symmetric and nonnegative. Also, for any $\alpha <1$%
, $\mu -\mu \in \mathcal{S}_{\alpha }$, so $D(\mu ,\mu )=0$. Now we only
need to check the triangle inequality:%
\begin{equation}
D(\mu _{1},\mu _{2})+D(\mu _{2},\mu _{3})\geq D(\mu _{1},\mu _{3}),
\label{triang}
\end{equation}%
for any probability measures $\mu _{1},\mu _{2},\mu _{3}$ on $\mathbb{R}^{m}$%
.

Note that (\ref{triang}) is clearly true if either $D(\mu _{1},\mu _{2})$ or
$D(\mu _{2},\mu _{3})$ is infinity. Thus we may assume that both $D(\mu
_{1},\mu _{2})$ and $D(\mu _{2},\mu _{3})$ are finite. So, $\mu _{1}-\mu
_{2}\in \mathcal{S}_{\alpha _{1}}$ and $\mu _{2}-\mu _{3}\in \mathcal{S}%
_{\alpha _{2}}$ for some $\alpha _{1},\alpha _{2}<1$. Now,
\begin{equation*}
\mu _{1}-\mu _{3}=\left( \mu _{1}-\mu _{2}\right) +\left( \mu _{2}-\mu
_{3}\right) \in \mathcal{S}_{\max \left\{ \alpha _{1},\alpha _{2}\right\}
}\left( \mathbb{R}^{m}\right) .
\end{equation*}%
Therefore,
\begin{equation*}
D(\mu _{1},\mu _{3})\leq \frac{1}{1-\max \left\{ \alpha _{1},\alpha
_{2}\right\} }\leq \frac{1}{1-\alpha _{1}}+\frac{1}{1-\alpha _{2}}.
\end{equation*}%
By taking an infimum over $\alpha _{1}$ and $\alpha _{2}$ we obtain the
triangular inequality (\ref{triang}).
\end{proof}

In general, $D$ is not necessarily a metric. Indeed, for any two atomic
probability measures $\mathbf{a,b}$ , we have $\mathbf{a-b\in }\mathcal{S}%
_{\alpha }\left( \mathbb{R}^{m}\right) $ for any $\alpha <1$. Thus, $D\left(
\mathbf{a,b}\right) =0$ while $\mathbf{a}$ and $\mathbf{b}$ are not
necessarily the same measure. Nevertheless, we may easily extend the
pseudometric $D$ to a metric on equivalent classes of measures. To this end
we define a notion of the equivalence classes on measures.

\begin{definition}
For any two probability measures $\mu $ and $\nu $ on $\mathbb{R}^{m}$, we
say
\begin{equation*}
\mu \sim \nu \text{\ if }D(\mu ,\nu )=0.
\end{equation*}%
The equivalent class of $\mu $ is denoted by $[\mu ]$.
\end{definition}

For instance, all atomic probability measures are equivalent to each other.

\begin{proposition}
If $\mu _{1}\sim \mu _{2}$ and $\nu _{1}\sim \nu _{2}$, then $D(\mu _{1},\nu
_{1})=D(\mu _{2},\nu _{2}).$
\end{proposition}

\begin{proof}
By definition of equivalence of measures we have $D(\mu _{1},\mu _{2})=0$
and $D(\nu _{1},\nu _{2})=0$. So by the triangular inequality (\ref{triang})
we have
\begin{equation*}
|D(\mu _{1},\nu _{1})-D(\mu _{2},\nu _{2})|\leq D(\mu _{1},\mu _{2})+D(\nu
_{1},\nu _{2})=0.
\end{equation*}
\end{proof}

%By this proposition, we have the following definition.

\begin{definition}
For any equivalent class $[\mu]$ and $[\nu]$, define
\begin{equation*}
\mathbf{D}([\mu],[\nu])=D(\mu,\nu).
\end{equation*}
\end{definition}

From this definition and proposition \ref{pseudometric}, clearly, we have
the following theorem.

\begin{theorem}
$\mathbf{D}$ is a metric on $\mathcal{P}\left( \mathbb{R}^{m}\right) /\sim $
.
\end{theorem}

\begin{definition}
The metric $\mathbf{D}$ is called the dimensional distance on the space $%
\mathcal{P}\left( \mathbb{R}^{m}\right) /\sim $ \ of equivalent classes of
probability measures on $\mathbb{R}^{m}$.
\end{definition}

We now give a geometric meaning to transport dimension of measures.

\begin{theorem}
For any positive probability measure $\mu $, we have
\begin{equation*}
\dim _{T}(\mu )=\mathbf{D}\left( \left[ \mu \right] ,\left[ \mathbf{a}\right]
\right) =D(\mu ,\mathbf{a})
\end{equation*}%
where $\mathbf{a}$ is any atomic probability measure.
\end{theorem}

\begin{proof}
Since $\mathbf{a\in }\mathcal{D}_{\alpha }$ for any $\alpha <1$, we have
\begin{eqnarray*}
D(\mu ,\mathbf{a}) &=&\inf_{\alpha <1}\{\frac{1}{1-\alpha }:\mu -\mathbf{a}%
\in S_{\alpha }\} \\
&=&\inf_{\alpha <1}\{\frac{1}{1-\alpha }:\mu \in \mathcal{D}_{\alpha
}\}=\dim _{T}\left( \mu \right) .
\end{eqnarray*}
\end{proof}

This theorem says that the transport dimension of a probability measure $\mu
$ is the distance from $\mu $ to any atomic measure with respect to the
dimensional distance. In other words, the dimension information of a measure
tells us quantitatively how far the measure is from being an atomic measure.

\end{document}